\documentclass[10pt]{article}

\usepackage[margin=2cm,bmargin=3cm]{geometry}
\usepackage{import, subcaption}
\usepackage{amsmath,amsfonts}
\usepackage{pgfplots}
\usepackage{color}
\usepackage{bm}
\usepackage{tikz}
\usepackage{authblk}
\usepackage{booktabs}
\usepackage{appendix}
\usepackage[backend=biber,bibencoding=utf8,maxcitenames=2,maxbibnames=9,style=numeric,
firstinits=true,url=false, isbn=false,date=year]{biblatex}
\AtBeginBibliography{\renewcommand*{\mkbibnamefamily}[1]{\textsc{#1}}}
\renewbibmacro*{doi+eprint+url}{%
	\printfield{doi}%
	\newunit\newblock%
	\iftoggle{bbx:eprint}{%
		\usebibmacro{eprint}%
	}{}%
	\newunit\newblock%
	\iffieldundef{doi}{%
		\usebibmacro{url+urldate}}%
	{}%
      }%
\renewcommand*{\bibfont}{\footnotesize}

\usepackage[hidelinks]{hyperref} 
\usepackage{doi}
\usetikzlibrary{calc}

\renewcommand\d[0]{\ensuremath{\operatorname{d}\!}}
\renewcommand\vec[0]{\bm}

\title{Finite volume discretization for poroelastic media with fractures modeled by
	contact mechanics}
\author[1]{Runar L. Berge} 
\author[1,2]{Inga Berre} 
\author[1]{Eirik Keilegavlen}
\author[1]{Jan M. Nordbotten}
\author[3]{Barbara Wohlmuth}

\affil[1]{Department of Mathematics, University of Bergen, Norway}
\affil[2]{NORCE, Norway}
\affil[3]{Department of Mathematics, Technische
		Universit\"at M\"unchen, Germany}

\date{}
\begin{document}
\maketitle
\begin{abstract}
  A fractured poroelastic body is considered where the opening of the fractures is
  governed by a nonpenetration law while slip is described by a Coulomb-type friction
  law. This physical model results in a nonlinear variational inequality problem. The
  variational inequality is rewritten as a complementary function, and a semismooth
  Newton method is used to solve the system of equations. For the discretization, we use
  a hybrid scheme where the displacements are given in terms of degrees of freedom per
  element, and an additional Lagrange multiplier representing the traction is added on
  the fracture faces. The novelty of our method comes from combining the Lagrange
  multiplier from the hybrid scheme with a finite volume discretization of the
  poroelastic Biot equation, which allows us to directly impose the inequality
  constraints on each subface.  The convergence of the method is studied for several
  challenging geometries in 2d and 3d, showing that the convergence rates of the finite
  volume scheme do not deteriorate when it is coupled to the Lagrange multipliers.  Our
  method is especially attractive for the poroelastic problem because it allows for a
  straightforward coupling between the matrix deformation, contact conditions, and fluid
  pressure.
\end{abstract}

\section{Introduction}
\label{sec:introduction}
Slip and opening of fractures due to fluid injection is of relevance to a number of
subsurface engineering processes. In hydraulic reservoir stimulation, the effect is
deliberately induced, while in storage operations and wastewater disposal, avoiding
reactivation and opening of fractures is important for preserving caprock integrity. In
any circumstance, triggering of larger slip events in the form of elevated levels of
seismicity must be avoided. The mathematical model of fracture resistance, slip and
opening results in a strongly coupled nonlinear variational inequality, which requires
advanced numerical schemes to solve. The purpose of this work is to describe and
implement a numerical method to solve this problem considering a poroelastic matrix.
The fractures are a set of predefined surfaces in the domain, and the nucleation or
growth of fractures is not considered.

The flow and mechanics of poroelastic media and the contact mechanics of elastic bodies
are well-developed research fields. For a porous or poroelastic medium, we refer to the
classical textbooks~\cite{coussy2003poromechanics,bear1988dynamics}. There exists an
extensive number of discretizations for the elliptic equations describing fluid flow in
a porous medium, and they all have different merits. The most popular discretizations
are the so-called locally conservative discretizations~\cite{russel2000relationships},
which include mixed finite elements~\cite{raviart1977mixed}, control-volume finite
elements~\cite{cai1997control-volume}, and finite volume
methods~\cite{aavatsmark2002introduction}. For the coupled poroelastic problem, without
considering fractures, it is known that a naive discretization of the coupling terms of
the fluid pressure and the solid displacement leads to stability issues for finite
element schemes~\cite{vermeer1981accuracy}. Several different methods have been proposed
to remove these
oscillations~\cite{haga2012on,murad1994stability,aguilar2008numerical}. Recently, a
finite volume method called the multipoint stress approximation was introduced for
elastic deformations~\cite{nordbotten2014cell-centered,keilegavlen2017finite}. This
method has been extended to the poroelastic Biot equations and shown to be stable
without adding any artificial stabilization terms in the limit of incompressible fluids
and small time steps~\cite{nordbotten2016stable}.

The contact mechanics problem, i.e., contact between two elastic bodies, is also the
topic of several textbooks~\cite{kikuchi1988contact,wriggers2006computational}. A widely
used solution strategy for the nonlinear variational inequalities resulting from the
mathematical formulation is the penalty method~\cite{kikuchi1981penalty}. The basic idea
is to penalize a violation of the inequality by adding extra energy to the system. The
solution depends then, in a very sensitive way, on the choice of the penalty
parameter. If the value of the parameter is too small, the condition number of the
algebraic system is extremely poor, and the nonlinear solver converges slowly. If the
value is too large, the accuracy of the solution is very poor, and unphysical
approximations can be obtained. Therefore, variationally consistent hybrid formulations
have gained interest recently.
The hybrid formulations are based on the discretization of additional unknown Lagrange
multipliers added to the contact region. This method has been applied to, among others,
the Signorini problem~\cite{belgacem2003hybrid}, frictional
contact~\cite{mcDevitt2000mortar}, and large deformations~\cite{puso2004finite}; see the
survey contribution~\cite{wohlmuth2011variationally} and the references therein.

For a poroelastic domain including fractures, different models for the contact problem
are
developed~\cite{moinfar2013coupled,mikelic2015phasefield,giovanardi2017unfitted,Garipov2016discrete}.
Most of these models, however, do not take into account the contact problem either by
assuming the fractures stick together~\cite{moinfar2013coupled} or that the fluid
pressure inside the fractures is so large that the fracture surfaces are never in
contact~\cite{mikelic2015phasefield,giovanardi2017unfitted}. The full contact problem for
a fractured poroelastic domain is considered by Garipov et al~\cite{Garipov2016discrete},
where they applied the penalty method to solve the nonlinear variational inequalities
resulting from the contact problem.

In the current work, we present a different numerical solution approach for poroelastic
media with contact mechanics. The discretization is based on a finite volume method for
poroelasticity~\cite{nordbotten2016stable} combined with a variationally consistent
hybrid discretization~\cite{hueber2008primal-dual,wohlmuth2011variationally}. The hybrid
formulation considered in this work can be regarded as a mortar
formulation~\cite{bernardi1994new} using matching meshes with the displacement as the
primal variable and the surface traction as the dual variable. The finite volume scheme
has previously been extended to fracture deformation by adding additional displacement
unknowns on the fracture faces~\cite{ucar2018finite}. This formulation was successfully
used to implement a fixed-point type iteration to approximate the friction
bound~\cite{ucar2017postinjection}; however, this formulation suffers from the fact that
a step length parameter needs to be tuned and that it might require many iterations to
converge~\cite{berge2019reactivation}.  An advantage of the scheme used in this work,
where the Lagrange multiplier of the hybrid formulation is coupled with the surface
traction obtained from the finite volume scheme, is that it gives a natural formulation
of the contact condition per subface. This formulation allows us to rapidly solve the
resulting nonlinear inequality problem by applying a semismooth Newton method; see the
work by H\"ueber et al~\cite{hueber2008primal-dual}, among
others~\cite{wohlmuth2011variationally,hueeber08discretization}.

The remainder of this paper is structured as follows. First, we state the problem and
give the governing equations. Then, the discretization is presented, which is divided
into two parts: (i) the finite volume discretization for the Biot equations and (ii) the
discrete hybrid formulation for the contact problem. We present four numerical
examples. The first two consider the dry case where the coupling between fluid pressure
and deformation of the rock is disregarded. The last two examples solves the poroelastic
deformation of a 2d and 3d domain where the deformation of fractures is governed by a
Coulomb friction law. Finally, we give concluding remarks.

\section{Problem statement}
\begin{figure}
  \centering \def\svgwidth{0.4\textwidth} 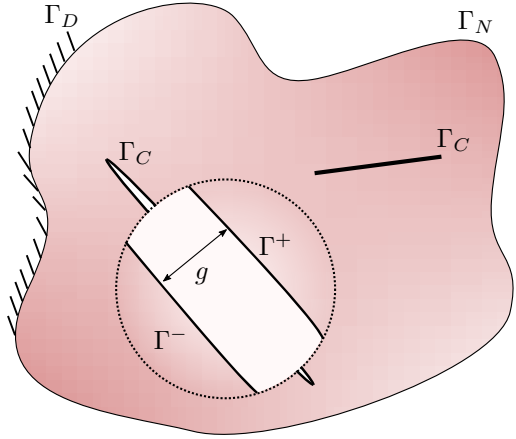
  \caption{A domain $\Omega$ where the external boundary is divided into two parts:
    $\Gamma_D$ and $\Gamma_N$. Included in the domain are two internal boundaries, or
    fractures, $\Gamma_C$. The two sides of the internal boundaries are labeled
    $\Gamma^+$ and $\Gamma^-$, as shown in the magnified circular region of the
    domain. The function $g(\vec x),\vec x \in\Gamma^+$ gives the initial gap between
    the two fracture sides. The left fracture has an initial gap $g>0$, while the top
    right fracture has an initial gap $g=0$. \label{fig:boundaries}}
	
\end{figure}

Let $\Omega$ be a fractured deformable porous body. The boundaries of the domain
$\partial \Omega$ are divided into three disjoint open sets,
$\Gamma_{D},\; \Gamma_{N},\; \textrm{and} \; \Gamma_{C}$, as illustrated in
Figure~\ref{fig:boundaries}: for the first set, a Dirichlet boundary condition is
assigned; for the second, a Neumann boundary condition is assigned; and the last is the
internal fracture boundary. We consider the Biot model for a poroelastic
medium~\cite{biot1941general}:
\begin{equation}\label{eq:biot}
  \begin{aligned}
    -\nabla \cdot \sigma&= \vec f_u &&\text{in } \Omega,\\
    \mathcal C : (\nabla \vec u + (\nabla \vec u)^\top) / 2 - \alpha p I &= \sigma
    &&\text{in } \Omega,\\
    c_0\dot p + \alpha\nabla\cdot \dot{\vec u} + \nabla\cdot \vec q &= f_p
    &&\text{in } \Omega,\\
    \vec q &= -\mathcal K\nabla p &&\text{in } \Omega,\\
    \vec u &= \vec g_{u,D} &&\text{on }\Gamma_{u,D},\\
    \sigma \cdot \vec n &= \vec g_{u,N} && \text{on } \Gamma_{u,N} ,\\
    p &= g_{p,D} &&\text{on }\Gamma_{p,D},\\
    \vec q \cdot \vec n &= g_{p,N} && \text{on } \Gamma_{p,N}.
  \end{aligned}
\end{equation}
The variables $\dot p$ and $\dot{\vec u}$ are the time derivatives of the pressure and
displacement, respectively. Throughout this paper we apply a backward Euler time
stepping, and $\dot p$ and $\dot{\vec u}$ should be interpreted as the discrete
derivatives
\begin{equation}\label{eq:backwar_euler}
  \dot p = \frac{p - p^{i}}{\Delta t},\quad \dot{\vec u} = \frac{\vec u - \vec
    u^{i}}{\Delta t},
\end{equation}
for the previous time iterate $i$ and time step length $\Delta t$. Note that we have
dropped the index for the current time iterate, $i+1$.  All parameters are, in general,
functions of space, e.g., $\mathcal C = \mathcal C(\vec x),\ \vec x\in \Omega$; however,
the explicit dependence is suppressed to keep the notation simple. Parameters associated
with the pressure $p$ and displacement $\vec u$ are given a subscript with the same
symbol. The vector $\vec f_u$ is a given body force, while $f_p$ is a given source
term. The stiffness tensor is denoted $\mathcal C$, the Biot coupling coefficient
$\alpha$, the storage coefficient $c_0$, and the permeability $\mathcal K$. Indicated by
subscripts, $\vec g$ represents Dirichlet and Neumann boundary conditions for
displacement and pressure, while $\vec n$ is the outwards pointing normal vector. In
this work, we use
$\mathcal C: (\nabla\vec u + (\nabla\vec u)^\top)/2 = G(\nabla\vec u + (\nabla\vec
u)^\top) + \Lambda\text{tr}(\nabla\vec u)I$, where $G$ and $\Lambda$ are the Lam\'e
parameters. Traction can also be derived for other material laws.

The fracture boundary, $\Gamma_C$, is divided into a positive side $\Gamma^+$ and a
negative side $\Gamma^+$. The choice of which side is positive and which is negative is
arbitrary and will only make a difference in the implementation. For the fracture
segments, a nonpenetration condition is enforced in the normal direction, meaning that
the positive and negative sides cannot penetrate each other. In the tangential
direction, a Coulomb friction law divides the contact region into a sliding part and a
sticking part. To formulate these contact conditions, the normal vector for the contact
region is defined as the normal vector of the positive side
$\vec n(\vec x) = \vec n^+(\vec x)$. Further, let
\begin{equation}\label{eq:positive_to_mortar}
  R : \Gamma^+ \rightarrow \Gamma^-
\end{equation}
be a mapping that projects a point from the positive boundary onto the negative boundary
as given by the normal vector. The gap function, which will appear in the nonpenetration
condition, is then defined as
\begin{equation*}
  g(\vec x) = \lVert\vec x - R(\vec x)\rVert \quad \vec x\in\Gamma^+,
\end{equation*}
where $\lVert\cdot\rVert$ is the Euclidean norm. Due to Newton's third law, the surface
traction, $\vec T = \sigma \cdot \vec n$, on the contact boundaries must be equal up to
the sign
\begin{equation}
  \label{eq:newtons_law}
  \vec T^+(\vec x) = -\vec T^-(R(\vec x)) \qquad \vec x\in \Gamma^+,
\end{equation}
and we use the notation $\vec T_C = \vec T^+$.  The surface traction is divided into a
normal and tangential part by
\begin{equation}\label{eq:normal_tangential}
  T_n(\vec x) = \vec T_C(\vec x)\cdot \vec n(\vec x),\quad
  \vec T_\tau(\vec x) = \vec T_C(\vec x) - T_n(\vec x)\vec n(\vec x) \qquad \vec x \in 
  \Gamma^+,
\end{equation}
and the displacement jump is defined as
$[\vec u(\vec x) ] = \vec u(\vec x) - \vec u(R(\vec x))$ for $\vec x\in \Gamma^+$. The
normal and tangential displacement jump is defined analogously to
Equation~\eqref{eq:normal_tangential}:
\begin{equation*} [\vec u(\vec x)]_n = [\vec u(\vec x)]\cdot \vec n(\vec x),\quad [\vec
  u(\vec x)]_\tau = [\vec u(\vec x)] - [\vec u(\vec x)]_n\vec n(\vec x) \qquad \vec x
  \in \Gamma^+.
\end{equation*}
The nonpenetration condition can now be formulated as
\begin{align}\label{eq:non_penetration}
  \begin{cases}
    [\vec u(\vec x)]_n - g(\vec x)  \leq 0 \\
    T_n(\vec x)([\vec u(\vec x)]_n - g(\vec x)) = 0 \\
    T_n(\vec x) \leq 0
  \end{cases}
  \qquad \vec x \in \Gamma^+,
\end{align}
where the first condition ensures that the two sides of the fracture cannot penetrate,
the second ensures that either the normal traction is zero or the fracture sides are in
contact, and the last enforces a negative normal component of the surface traction.

The tangential part of the surface traction is governed by a Coulomb friction law:
\begin{align}\label{eq:coulomb_friction}
  \begin{cases}
    \lVert\vec T_\tau(\vec x)\rVert &\leq F(\vec x)|T_n(\vec x)| \\
    \lVert\vec T_\tau(\vec x)\rVert & < F(\vec x)|T_n(\vec x)|\ \rightarrow [\dot{\vec
      u}(\vec x)]_\tau = 0 \\
    \lVert\vec T_\tau(\vec x) \rVert &= F(\vec x)|T_n(\vec x)|\ \rightarrow \exists\
    \zeta \in \mathbb R : \vec T_\tau(\vec x) = -\zeta^2[\dot{\vec u}(\vec x)]_\tau
  \end{cases}
                                       \qquad \vec x \in \Gamma^+,
\end{align}
where $F$ is the coefficient of friction, and $\dot{\vec u}$ is the displacement
velocity approximated by the backward Euler scheme, as given by
Equation~\eqref{eq:backwar_euler}. The first equation gives the friction bound, the
second ensures that if the friction bound is not reached, then the surface is sticking,
and the last equation ensures that if the friction bound is reached, then the tangential
sliding velocity is parallel to the tangential traction. In the static case, e.g., for
the purely mechanical problem when $\alpha=0$, the notion of a velocity does not
exist. For these cases, it is common to replace the sliding velocity,
$[\dot{\vec u}]_\tau$, by the displacement jump, $[\vec u]_\tau$, in
Equation~\eqref{eq:coulomb_friction}~\cite{wohlmuth2011variationally}.

For the fluid, the fractures are modeled as impermeable. This means that the fluid
cannot flow in or across the fractures, i.e.,
$\vec q(\vec x)\cdot \vec n(\vec x) = 0,\ \vec x \in \Gamma_C$. To avoid excessive model
complexity, we have chosen a model with impermeable fractures. For possible methods to
extend this work to include fracture flow, we refer to the work by Dietrich et
al~\cite{dietrich2005}, among
others~\cite{ucar2018,Garipov2016discrete,martin2005,angelo2012}.

\section{Discretization}\label{sec:discretization}
We define the triplet $(\mathcal T, \mathcal F, \mathcal V)$ as the cells, faces and
vertices of our mesh. It is assumed that the mesh conforms to the fractures; that is,
the positive and negative sides of the fractures are tessellated by a subset of the
faces of the mesh. Before the discretization is described, we need to define some
notation, and we start by giving the relation between cells, faces and vertices using
the standard notation for finite-volume
methods~\cite{eymard2000,nordbotten2016stable}:
\begin{itemize}
\item For a cell $K\in \mathcal T$, we denote its faces by $\mathcal F_K$ and its
  vertices by ~$\mathcal V_K$.
\item For a face $\pi\in \mathcal F$, we denote the neighboring cells as
  $\mathcal T_\pi$ and its vertices as $\mathcal V_\pi$.
\item For a vertex $v\in\mathcal V$, we denote the adjacent cells by $\mathcal T_v$ and
  the adjacent faces by $\mathcal F_v$.
\end{itemize}
In addition to the mesh triplet $(\mathcal T, \mathcal F, \mathcal V)$, we define the
so-called subcells and subfaces illustrated in Figure~\ref{fig:sub_cell}:
\begin{figure}
  \centering \def\svgwidth{0.3\textwidth}
\begingroup%
  \makeatletter%
  \providecommand\color[2][]{%
    \errmessage{(Inkscape) Color is used for the text in Inkscape, but the package 'color.sty' is not loaded}%
    \renewcommand\color[2][]{}%
  }%
  \providecommand\transparent[1]{%
    \errmessage{(Inkscape) Transparency is used (non-zero) for the text in Inkscape, but the package 'transparent.sty' is not loaded}%
    \renewcommand\transparent[1]{}%
  }%
  \providecommand\rotatebox[2]{#2}%
  \newcommand*\fsize{\dimexpr\f@size pt\relax}%
  \newcommand*\lineheight[1]{\fontsize{\fsize}{#1\fsize}\selectfont}%
  \ifx\svgwidth\undefined%
    \setlength{\unitlength}{403.39794249bp}%
    \ifx\svgscale\undefined%
      \relax%
    \else%
      \setlength{\unitlength}{\unitlength * \real{\svgscale}}%
    \fi%
  \else%
    \setlength{\unitlength}{\svgwidth}%
  \fi%
  \global\let\svgwidth\undefined%
  \global\let\svgscale\undefined%
  \makeatother%
  \begin{picture}(1,0.82538911)%
    \lineheight{1}%
    \setlength\tabcolsep{0pt}%
    \put(0,0){\includegraphics[width=\unitlength,page=1]{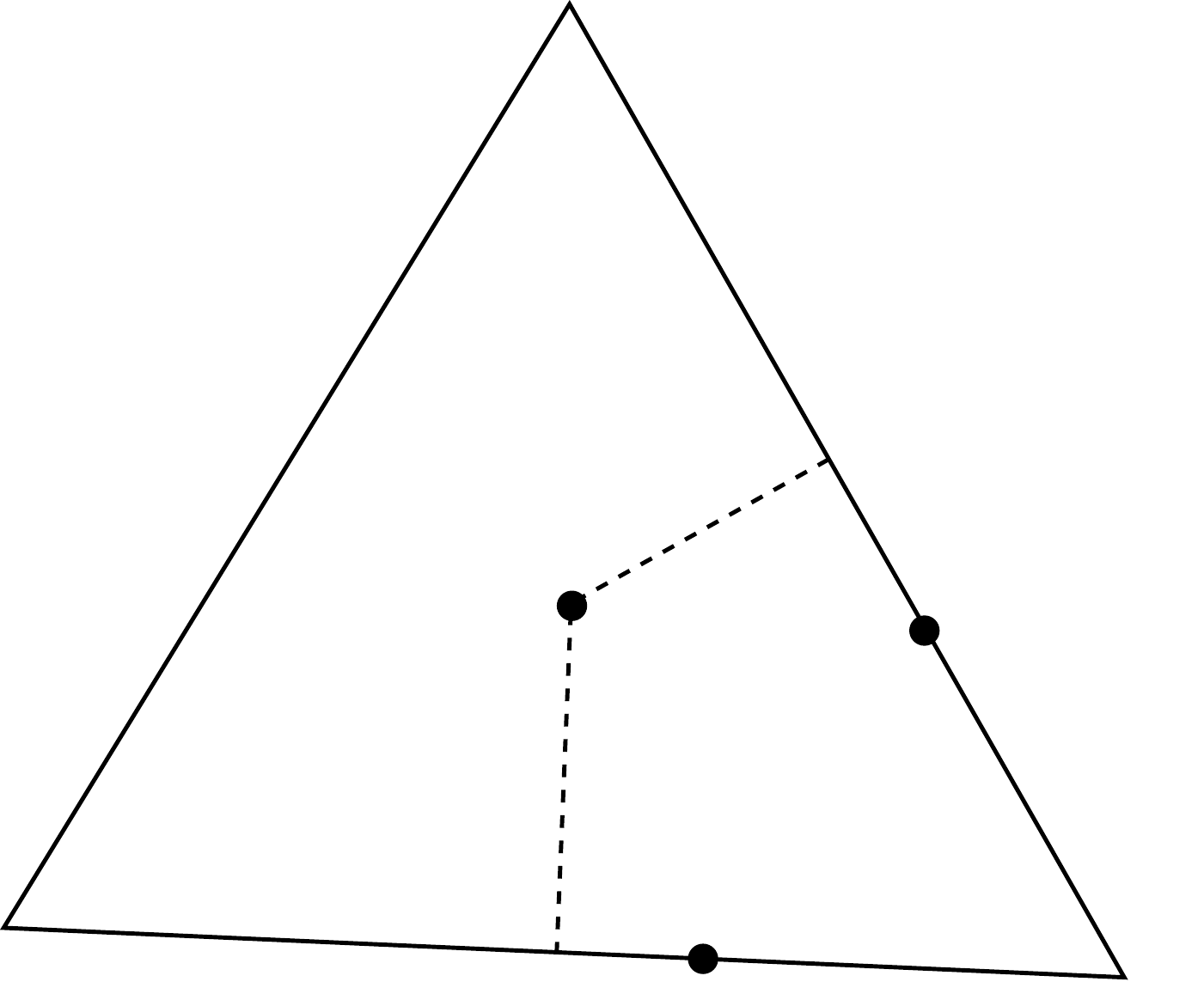}}%
    \put(0.45118836,0.31995073){\color[rgb]{0,0,0}\makebox(0,0)[rt]{\lineheight{1.25}\smash{\begin{tabular}[t]{r}$\vec x_K$\end{tabular}}}}%
    \put(0.69230068,0.47341652){\color[rgb]{0,0,0}\makebox(0,0)[lt]{\lineheight{1.25}\smash{\begin{tabular}[t]{l}$\pi$\end{tabular}}}}%
    \put(0.44489366,0.44915769){\color[rgb]{0,0,0}\makebox(0,0)[lt]{\lineheight{1.25}\smash{\begin{tabular}[t]{l}$K$\end{tabular}}}}%
    \put(0,0){\includegraphics[width=\unitlength,page=2]{sub_cell.pdf}}%
    \put(0.96037297,0.00993671){\color[rgb]{0,0,0}\makebox(0,0)[lt]{\lineheight{1.25}\smash{\begin{tabular}[t]{l}$v$\end{tabular}}}}%
    \put(0.86608467,0.20382553){\color[rgb]{0,0,0}\makebox(0,0)[lt]{\lineheight{1.25}\smash{\begin{tabular}[t]{l}$\pi_v$\end{tabular}}}}%
    \put(0.79171645,0.33529799){\color[rgb]{0,0,0}\makebox(0,0)[lt]{\lineheight{1.25}\smash{\begin{tabular}[t]{l}$\vec x_\pi^v$\end{tabular}}}}%
    \put(0.63501233,0.17885828){\color[rgb]{0,0,0}\makebox(0,0)[lt]{\lineheight{1.25}\smash{\begin{tabular}[t]{l}$K_v$\end{tabular}}}}%
  \end{picture}%
\endgroup%

  \caption{Notation used to describe the mesh. For a cell $K$, face $\pi$ and vertex $v$
    of the mesh, we associate a subcell $(K,v)$ and subface $(\pi,v)$, as well as a cell
    center $\vec x_K$ and continuity point $\vec x_\pi^v$. In this figure, the cell is
    the full triangle, and the subcell is given by the gray area. \label{fig:sub_cell}}
\end{figure}
\begin{itemize}
\item For a vertex $v\in \mathcal V_K$, we define a subcell of $K$ identified by
  $(K, v)$ with a volume $m_K^v$ such that
  $\sum_{v\in\mathcal V_K} m_K^v = m_K=\int_K\d \vec x$.
\item For a vertex $v\in \mathcal V_\pi$, we associate a subface identified by
  $(\pi, v)$ with an area $m_\pi^v$ such that
  $\sum_{v\in\mathcal V_\pi}m_\pi^v = m_\pi=\int_\pi\d\vec x$.
\end{itemize}
The subfaces cannot be chosen arbitrary but should correspond to faces of the subcells;
for the triplet $(v, \pi, K)$, the intersection of the boundary of the subcell $(K, v)$
and the face $\pi$ should equal the subface $(\pi, v) = \partial (K, v)\cap
\pi$. Further, all subcells and subfaces are assumed to have a positive measure. Note
that in an abuse of notation, we use $K$ for both indexing and the geometric object so
that both $\mathcal V_K$ and $\int_K\d \vec x$ make sense. All subfaces
$(\pi, v),\ \pi \in \mathcal F,\ v\in\mathcal V_\pi$ are divided into three disjoint
sets $\mathcal P,\mathcal N$, and $\mathcal R$, where $\mathcal P$ contains all subfaces
located on the positive boundary $\Gamma^+$, $\mathcal N$ contains all subfaces located
on the negative boundary $\Gamma^-$, and $\mathcal R$ contains the remaining subfaces.

Finally, for each element $K\in \mathcal T$, a cell center $\vec x_K\in K$ is defined,
and for each subface $(\pi, v)$, we associate a continuity point $\vec x_\pi^v$ located
at any point on the subface, $(\pi, v)$, however, the distance to the vertex $v$ must be
greater than zero. The unit normal for each face is denoted by $\vec n_\pi$, which is
equal to the subface normal of the face $\vec n_\pi^v$. When it is necessary to
distinguish the direction of the normal, it is defined as the outward pointing normal
$\vec n_K^\pi$ of a cell $K\in\mathcal T_\pi$. Note that for a face $\pi$, we have
$\mathcal T_\pi = \{K, L\}$, $\vec n_K^\pi = -\vec n_L^\pi$.

In the implementation used in the examples of this paper, the following construction is
employed: The face- and cell-centers are chosen as the centroid of the corresponding
face and cell. In 2d, the subface $(\pi, v)$ is defined by the convex-hull of the vertex
$v$ and the face-center $\vec x_\pi$. In 3d, the subface is in addition defined by the
midpoints of the edges of the face $\pi$ that are connected to $v$. For simplices, this
construction partition each face into a set of subfaces of equal area. The subcell
$(K, v)$ is defined by the convex-hull of the cell-center $\vec x_K$ and the subfaces
$(\pi, v),\ \pi\in \mathcal F_v\cap \mathcal F_K$. The continuity point, $\vec x_\pi^v$,
is taken to be one third the distance from the face-center to the vertex,
$\vec x_\pi^v = \vec x_\pi - (\vec x_\pi - v) / 3$. An example of this construction is
shown in Figure~\ref{fig:sub_cell}.

\subsection{Finite volume discretization}
We use a finite volume discretization~\cite{nordbotten2016stable} to discretize the Biot
Equations~\eqref{eq:biot}.  This is based on two discrete variables, $\vec u_K$ and
$p_K$, which are the cell-centered displacement and pressure, respectively. Within each
subcell $(K,v),\ K\in \mathcal T,\ v \in \mathcal V_K$, it is assumed that the
displacements and pressures are linear in each subcell, and the discrete gradients are
denoted by $(\bar\nabla \vec u)_K^v$ and $(\bar\nabla p)_K^v$, where the bar over the
gradient operator is added to distinguish it from the continuous gradients. For the
mechanical stress, we adapt the notion of weak symmetry~\cite{keilegavlen2017finite};
given the volume weighted average
\begin{equation*}
  \left<\Xi\right>_v =\frac{1}{ \sum_{K\in \mathcal T_v}m_K^v} \sum_{K\in \mathcal
    T_v}m_K^v\Xi_K^v,
\end{equation*}
associated with a vertex $v$, the discrete weakly symmetric mechanical stress is given
by
\begin{equation}\label{eq:weak_symmetry}
  \vec \theta_K^v = \mathcal C_K:(\bar\nabla \vec u)_K^v - 
  \frac{\left<\mathcal C : (\bar\nabla \vec u)\right>_v 
    - \left<\mathcal C : (\bar\nabla \vec u)\right>_v^\top}{2}.
\end{equation}
This is referred to as weak symmetry because
\begin{equation*}
  \left< \vec \theta - \vec \theta^\top\right>_v = 0.
\end{equation*}
To simplify notation, the tensor $\mathcal C_K^v$ is referred to as the stress tensor,
which acts to weakly symmetrize the stress:
\begin{equation*}
  \vec \theta_K^v = C_K^v:(\bar\nabla\vec u)_K^v.
\end{equation*}
The expression $C_K^v:(\bar\nabla\vec u)_K^v$ should not be interpreted as a single
tensor vector product but as a weighted sum of products given by
Equation~\eqref{eq:weak_symmetry}.

Using the weak symmetry, the flux and traction over each subface given by the discrete
variables can be stated as
\begin{align}\label{eq:darcy_flux}
  q_{K,\pi}^v &= -m_\pi^v\mathcal K_K(\bar\nabla p)_K^v \cdot \vec n_K^\pi,\\
  \label{eq:biot_stress}
  \vec T_{K,\pi}^v &= m_\pi^v(\mathcal C_K^v:(\bar \nabla \vec u)_K^v - \alpha p_K I)
                     \cdot \vec n_K^\pi.
\end{align}
For a spatially varying permeability and stress tensor, we use the cell-center value to
evaluate the parameters $\mathcal K_K = \mathcal K(\vec x_K)$ for each cell.

The finite-volume scheme will be constructed such that the gradient unknowns can be
eliminated by performing a local static condensation. The following presentation on how
we obtain the numerical gradients is designed to reflect the computer
implementation. This presentation is different from most classical papers on the used
finite-volume discretization~\cite{aavatsmark2002introduction}, however, it is
mathematically equivalent.  After the numerical gradients are expressed in terms of the
cell-center unknowns, the finite volume structure of the discretization is obtained by
enforcing mass/momentum conservation for each cell. The final scheme will be locally
conservative and given by the cell-centered displacement and pressure. A detail that
will be important when we introduce the hybrid discretization is the possibility of
exactly reconstructing the discrete gradients, and thus also the flux and traction, from
the cell-centered variables $\vec u$ and $p$.

The discrete fluid flux given in Equation~\eqref{eq:darcy_flux} does not contain any
dependence on the displacement, and it is identical to the fluid flux for the uncoupled
fluid pressure, i.e., $\alpha=0$. To discretize the flux, we use the MPFA-O scheme for
simplices~\cite{klausen2008,friis2009}. Each subcell gradient $(\bar\nabla p)_K^v$ is
associated with a fluid flux as given in~\eqref{eq:darcy_flux}. Conservation of mass is
enforced for each internal subface. This requires the fluid flux for cells
$(K, L) \in \mathcal T_\pi$ sharing a face $\pi$ to be equal and opposite over each of
their shared subfaces; that is,
\begin{equation}\label{eq:flux_cont}
  -m_\pi^v\mathcal K_K(\bar\nabla p)_K^v \cdot \vec n_K^\pi
  = m_\pi^v\mathcal K_L(\bar\nabla p)_L^v \cdot \vec n_L^\pi.
\end{equation}
The pressure is not required to be continuous across the whole subface. Instead,
pressure continuity is enforced at the continuity points, $\vec x_\pi^v$, that is,
\begin{equation}\label{eq:pressure_cont}
  p_K + (\bar\nabla p)_K^v \cdot (\vec x_\pi^v - \vec x_K) =
  p_L + (\bar\nabla p)_L^v \cdot (\vec x_\pi^v - \vec x_L).
\end{equation}
Here, we have made use of the assumption that the pressure is linear in each subcell to
write the pressure at the continuity point $\vec x_\pi^v$ as a function of the cell
center pressure $p_K$ and gradient $(\bar\nabla p)_K^v$. If a subface is on the Neumann
boundary, $\pi \subset \Gamma_{p,N}$, the flux over the subface is given by evaluating
the boundary condition at the continuity point and multiply the value with the subface
area:
\begin{equation}\label{eq:flux_bc}
  -m_\pi^v\mathcal K_K(\bar\nabla p)_K^v \cdot \vec n_K^\pi = m_\pi^vg_{p,N}(\vec x_\pi^v).
\end{equation}
If a subface is on the Dirichlet boundary, $\pi \subset \Gamma_{p,D}$, the pressure on
the subface is given by
\begin{equation}\label{eq:pressure_bc}
  p_K + (\bar\nabla p)_K^v \cdot (\vec x_\pi^v - \vec x_K) = g_{p,D}(\vec x_\pi^v).
\end{equation}
Faces on the fracture boundary, $(\pi, v)\subset \Gamma^+\cup\Gamma^-$, is given a zero
Neumann boundary condition as we have assumed impermeable fractures.

Around each vertex $v$ we can now form a local linear system of equations from which the
gradients $(\bar\nabla p)_K^v,\ K\in\mathcal T_v$ can be eliminated:
\begin{equation}\label{eq:gradp}
  \vec{(\bar\nabla p)}_v
  =
  \begin{bmatrix}
    Q_p\\
    D_{p,G}
  \end{bmatrix}^{-1} \left(
    \begin{bmatrix}
      \vec g_{p,N} \\\vec g_{p,D}
    \end{bmatrix}
    -
    \begin{bmatrix}
      0\\
      D_p
    \end{bmatrix}
    \vec p
  \right).
\end{equation}
The first block $Q_p\vec{(\bar\nabla p)}_v = \vec g_{p,N}$ in this linear system is the
collection of all flux balance Equations~\eqref{eq:flux_cont} and~\eqref{eq:flux_bc} for
the vertex $v$. The next block
$D_{p,G}\vec{(\bar\nabla p)}_v = \vec g_{p,D} - D_p\vec p$ collects all the pressure
continuity Equations~\eqref{eq:pressure_cont} and~\eqref{eq:pressure_bc}. Thus,
$\vec{ (\bar\nabla p)}_v$ is the vector of the subcell gradients $(\bar\nabla p)_K^v$,
the matrix $Q_p$ represents products of the form $m_\pi^v\vec n_K^\pi\mathcal K_K$, the
matrix $ D_{p,G}$ represents the distances $\vec x_\pi^v - \vec x_K$, the vectors
$\vec g_{p,N}$ and $\vec g_{p,D}$ are possible boundary conditions, and $ D_p$ has
entries $1$ for $p_K^v$ and $-1$ for $p_L^v$.

The elimination of the displacement gradients $(\bar\nabla\vec u)_K^v$ is similar to the
elimination of the pressure gradients $\bar \nabla p_K^v$. First, the continuity of
traction gives us for each internal subface
\begin{equation}\label{eq:stress_cont}
  m_\pi^v(\mathcal C_K^v:(\bar \nabla \vec u)_K^v - \alpha p_K I) \cdot \vec
  n_K^\pi =
  -m_\pi^v(\mathcal C_L^v:(\bar \nabla \vec u)_L^v - \alpha p_L I) \cdot \vec n_L^\pi.
\end{equation}
It is worth pointing out that, for internal faces, the averaging part of the operator
$\mathcal C_K^v : (\bar\nabla\vec u)_K^v$ is the same on the right- and left-hand sides.
Thus, the balance of traction can be written as
\begin{equation*}
  m_\pi^v(\mathcal C_K:(\bar \nabla \vec u)_K^v - \alpha p_K I) \cdot \vec
  n_K^\pi =
  -m_\pi^v(\mathcal C_L:(\bar \nabla \vec u)_L^v - \alpha p_L I) \cdot \vec n_L^\pi.
\end{equation*}
However, for boundary faces, the complete Equation~\eqref{eq:stress_cont} must be used.
Unlike the fluid fluxes in~\eqref{eq:darcy_flux}, the traction is different from the
uncoupled system due to the term $\alpha p_K I$. It is important to include the Biot
stress in the local systems to obtain the correct force balance in our
method~\cite{nordbotten2016stable}. We will see later that this approach also gives a
higher-order term in the mass balance for the fluid, which acts analogously to the
stabilization terms in other colocated schemes. For the fluid pressure, displacement
continuity is enforced at the continuity points $\vec x_\pi^v$:
\begin{equation}\label{eq:disp_cont}
  \vec u_K + (\bar \nabla \vec u)_K^v(\vec x_\pi^v - \vec x_K) =
  \vec u_L + (\bar \nabla \vec u)_L^v(\vec x_\pi^v - \vec x_L).
\end{equation}
For a subface on the Neumann boundary, $(\pi, v)\subset \Gamma_{u,N}$, the boundary
condition is evaluated at the continuity point and multiplied with the subface area:
\begin{equation}\label{eq:stress_bc}
  m_\pi^v(\mathcal C_K^v:(\bar \nabla \vec u)_K^v - \alpha p_K I) \cdot \vec
  n_K^\pi = m_\pi^v\vec g_{u,N}(\vec x_\pi^v),
\end{equation}
For a subface on the Dirichlet boundary $(\pi, v)\subset \Gamma_{u,D}$ the displacement
vector at the continuity point is given:
\begin{equation}\label{eq:disp_bc}
  \vec u_K + (\bar \nabla \vec u)_K^v(\vec x_\pi^v - \vec x_K) = \vec g_{u,D}(\vec x_\pi^v).
\end{equation}
Subfaces on the fracture boundary is given a Neumann condition defined by the Lagrange
multiplier:
\begin{equation}\label{eq:traction_contact}
  \begin{aligned}
    m_\pi^v(\mathcal C_K^v:(\bar \nabla \vec u)_K^v - \alpha p_K I) \cdot \vec
    n_K^\pi &= m_\pi^v\vec\lambda_\pi^v &&(\pi, v) \in\mathcal P,\\
    m_\pi^v(\mathcal C_K^v:(\bar \nabla \vec u)_K^v - \alpha p_K I) \cdot \vec n_K^\pi
    &= -m_\pi^v\vec\lambda(R^{-1}(\vec x_\pi^v)) &&(\pi, v) \in \mathcal N.
  \end{aligned}
\end{equation}
The contribution to the negative side $-\vec\lambda(R^{-1}(\vec x_\pi^v))$ is just the
mapping onto the Lagrange multiplier on the corresponding positive subface as given by
Equation~\eqref{eq:positive_to_mortar}.

A local elimination of the displacement gradients $(\bar\nabla\vec u)_K^v$ can now be
done around each vertex to express them in terms of the cell-center displacement and
pressure:
\begin{equation}\label{eq:gradu}
  \vec{ (\bar\nabla\vec u)_v}
  =
  \begin{bmatrix}
    Q_u\\
    D_{u,G}
  \end{bmatrix}^{-1} \left(
    \begin{bmatrix}
      \vec g_{u,N} \\%
      \vec g_{u,D}
    \end{bmatrix}
    -
    \begin{bmatrix}
      P & 0 & -M_\pm\\
      0 & D_U & 0
    \end{bmatrix}
    \begin{bmatrix}
      \vec p\\
      \vec u\\
      \vec\lambda
    \end{bmatrix}
  \right),
\end{equation}
The variable $\vec{(\bar\nabla u)_v}$ is the vector of the displacement gradients,
$(\bar\nabla\vec u)_K^v$, around the vertex $v$, the matrix $ Q_u$ represents products
of the form $m_\pi^v\vec n^\top\mathcal C_K^v$, the matrix $ D_{u, G}$ represents the
same distance vectors as in~\eqref{eq:gradp}, the vectors $\vec g_{u,N}$ and
$\vec g_{u,D}$ are possible boundary conditions, and $ D_u$ is a matrix with entries
$\pm 1$. The term $P$ is the only difference between the coupled and uncoupled system
and contains products of the form $m_\pi^v\alpha I\vec n_\pi^v$, and the matrix $M_\pm$
contains the positive areas $m_\pi^v$ for the positive subfaces and the negative areas
$-m_\pi^v$ for the negative subfaces and represents the Lagrange multiplier contribution
to the traction balance in Equation~\eqref{eq:traction_contact}.

The finite-volume discretization of fluid flow is then obtained by expressing the fluid
mass conservation over each cell~$K$ in terms of the discrete variables,
\begin{equation}\label{eq:flux_balance_cell}
  \sum_{\pi\in \mathcal F_K}\sum_{v\in \mathcal V_\pi} -m_\pi^v\mathcal 
  K_K(\bar\nabla p)_K^v \cdot \vec n_K^\pi + \sum_{v\in\mathcal V_K} [m_K^v\alpha 
  (\bar\nabla\cdot\dot{\vec u})_K^v
  +c_0 m_K^v\dot p_K] =
  \int_Kf_p\d\vec x.
\end{equation}
The pressure gradient $(\bar\nabla p)_K^v$ and displacement divergence
$(\bar\nabla\cdot\vec u)_K^v = \text{tr}(\bar\nabla\vec u)_K^v$ are obtained as linear
functions of the cell-centered pressures and displacements and Lagrange multipliers from
the local systems given in~\eqref{eq:gradp} and~\eqref{eq:gradu}. The appearance of the
pressure in the discrete displacement divergence is essential for the consistency of the
method and is similar to the artificially introduced stability terms in other methods;
see, e.g., Gaspar et al~\cite{gaspar2008stabilized}.

For the mechanics, momentum is conserved for all cells $K$,
\begin{equation}\label{eq:momentum_balance_cell}
  -\sum_{\pi\in\mathcal F_K}\sum_{v\in\mathcal V_\pi} m_\pi^v\mathcal 
  C_K^v:(\bar\nabla\vec
  u)_K^v\cdot \vec n_K^\pi = \int_K \vec f \d V.
\end{equation}
Note that the term $\alpha p_KI$ from the Biot stress in~\eqref{eq:biot_stress} sums to
zero over a cell due to Gauss's theorem; however, the pressure dependence on the subcell
gradients gives the correct fluid pressure contribution to the mechanics. Similarly, the
dependence of the Lagrange multiplier on the subcell gradients gives the correct force
contribution to the momentum balance.

To summarize, the finite volume scheme is constructed by defining a set of discrete
pressure and displacement gradients for each subcell. Flux and pressure continuity is
enforced over each subface for the fluid, and traction and displacement continuity is
enforced for each subface for the solid. This defines a small local system around each
node from which the pressure and displacement gradients can be expressed as a linear
combination of the cell-centered pressure and displacement, and Lagrange multiplier and
then eliminated. A stable coupling between the fluid and solid is achieved by
considering the Biot stress, i.e.,
$\mathcal C : (\nabla \vec u + (\nabla \vec u)^\top) / 2 - \alpha p I$, for traction
balance of the local systems.

\subsection{Hybrid formulation}\label{sec:hybrid_formulation}
To solve the contact conditions~\eqref{eq:non_penetration} and
\eqref{eq:coulomb_friction}, we apply the active-set strategy, which is equivalent to a
semismooth Newton method described by H\"ueber et al~\cite{hueber2008primal-dual}. See
also the paper by Wohlmuth~\cite{wohlmuth2011variationally}.  We recapitulate the
solution strategy in this section for the completeness of this paper. The main
difference in our approach is how the Lagrange multipliers, which represent the surface
traction, are coupled to the displacement unknowns in the surrounding domain. In our
finite volume scheme, the Lagrange multipliers enter into the local equations for the
displacement gradients.

A set of Lagrange multipliers is defined on the positive subface boundaries
\begin{equation*}
  \vec\lambda_\pi^v = \sigma_\pi^v\cdot \vec n_\pi^v, \qquad(\pi, v)\in\mathcal P.
\end{equation*}
The normal $\lambda_{\pi n}^v$ and tangential $\vec \lambda_{\pi\tau}^v$ components of
the Lagrange multiplier are defined analogously to~\eqref{eq:normal_tangential}. The
displacement on the subfaces, $\vec u_\pi^v$, is obtained as in
Equation~\eqref{eq:disp_bc} for local systems.

The discrete formulation of the nonpenetration condition~\eqref{eq:non_penetration} can
for each subface be written as
\begin{align}\label{eq:discrete_non_penetration}
  \begin{cases}
    [\vec u_\pi^v]_n - g_\pi^v  \leq 0 \\
    \lambda_{\pi n}^v([\vec u_\pi]_n^v - g_\pi^v) = 0 \\
    \lambda_{\pi n}^v \leq 0
  \end{cases}
  \qquad (\pi, v) \in \mathcal P,
\end{align}
and the static Coulomb friction~\eqref{eq:coulomb_friction} as
\begin{align}\label{eq:discrete_coulomb_friction}
  \begin{cases}
    \lVert\vec\lambda_{\pi\tau}^v\rVert &\leq F_\pi^v|\lambda_{\pi n}^v| \\
    \lVert\vec\lambda_{\pi\tau}^v\rVert & < F_\pi^v|\lambda_{\pi n}^v|\
    \rightarrow [\dot{\vec u}_\pi^v]_\tau = 0 \\
    \lVert\vec\lambda_{\pi\tau}^v \rVert &= F_\pi^v|\lambda_{\pi n}^v|\ \rightarrow
    \exists \zeta \in \mathbb R : \vec\lambda_{\pi\tau}^v = -\zeta^2[\dot{\vec
      u}_\pi^v]_\tau
  \end{cases}
                                           \qquad (\pi, v)\in\mathcal P.
\end{align}
Recall that for the static case, the sliding velocity is replaced by the displacement
jump, $[\vec u_\pi^v]_\tau$. We define
$b_\pi^{v,k} = F_\pi^v(-\lambda_{\pi n}^{v,k} + c([\vec u_\pi^{v,k}]_n - g_\pi^v))$,
which can be interpreted as the friction bound. The nonpenetration condition can now be
rewritten as the nonlinear complementary function
\begin{equation}\label{eq:comp_n}
  C_n([\vec u_\pi^v]_n, \lambda_{\pi n}^v)
  = -\lambda_{\pi n}^v -\frac{1}{F_\pi^v} \max\{0, b_\pi^{v,k}\},
\end{equation}
where $c>0$ is a given numerical parameter and $\max\{\cdot,\cdot\}$ is the maximum
function.  Similarly, we can now rewrite the Coulomb friction as the complementary
function
\begin{equation}\label{eq:comp_t}
  \begin{aligned}
    \vec C_\tau([\dot{\vec u}_\pi^v], \vec\lambda_{\pi\tau}^v) = \max\{b_\pi^{v,k},
    \lVert -\vec\lambda_{\pi\tau}^v + c [\vec u_\pi^v]_\tau\rVert
    \}(-\vec\lambda_{\pi\tau}^v) - \max\{0, b_\pi^{v,k}\}(-\vec\lambda_{\pi\tau}^v +
    c[\dot{\vec u}_\pi^v]_\tau).
  \end{aligned}
\end{equation}
The solution pair $(\vec u_\pi^v, \vec\lambda_\pi^v)$ satisfies the nonpenetrating
condition~\eqref{eq:discrete_non_penetration} and Coulomb
law~\eqref{eq:discrete_coulomb_friction} if and only if
$C_n([\vec u_\pi^v]_n , \lambda_{\pi n}^v) =0$ and
$\vec C_\tau([\dot{\vec u}_\pi^v], \vec\lambda_\pi^v)= \vec 0$. We apply a semismooth
Newton method to this problem, which results in an active set method. Given the solution
$(\vec u^k, \vec \lambda^k)$ from the previous Newton iteration, we divide the contact
subfaces into three disjoint sets:
\begin{equation}\label{eq:sliding_sets}
  \begin{aligned}
    \mathcal I^{k+1}_n &= \{(\pi, v)\in \mathcal P: b_\pi^{v,k} \le 0 \}\\
    \mathcal I^{k+1}_\tau &= \{(\pi, v)\in \mathcal P: \lVert
    -\vec\lambda_{\pi\tau}^{v,k} + c\vec
    [\dot{\vec u}_{\pi}^{v,k}]_\tau\rVert - b_\pi^{v,k} < 0\}\\
    \mathcal A^{k+1} &= \{(\pi, v)\in \mathcal P: \lVert -\vec\lambda_{\pi\tau}^{v,k} +
    c\vec [\dot{\vec u}_{\pi}^{v,k}]_\tau\rVert \ge b_\pi^{v,k} > 0\}.
  \end{aligned}
\end{equation}
The first set contains the subfaces not in contact. The second set contains the subfaces
in contact whose friction bound is not reached, i.e., they are sticking. The third set
contains the subfaces in contact where the friction bound is reached, i.e., they are
sliding. The new iterates $([\vec u_{\pi}^{v,k+1}]_\tau, \vec\lambda_{\pi\tau}^{v,k+1})$
in the semismooth Newton scheme are then calculated depending on which set the subface
belongs to. The update is found by calculating the derivative of the complementary
functions $C_n$ and $\vec C_\tau$ for each of the three sets. For the subfaces not in
contact, zero traction is enforced
\begin{equation}\label{eq:open}
  \vec\lambda_\pi^{v,k+1} = \vec 0 \qquad   (\pi, v) \in \mathcal I_n^{k+1}.
\end{equation}
For the subfaces in contact and sticking, we enforce
\begin{equation}\label{eq:stick}
  [u_{\pi n}^{v,k+1}] = g_\pi^v, \qquad 
  [\dot{\vec u}_{\pi}^{v,k+1}]_\tau + \frac{F_\pi^v
    [\dot{\vec u}_{\pi}^{v,k}]_\tau}{ b_\pi^{v,k}} \lambda_{\pi n}^{v,k+1} = [\dot{\vec
    u}_{\pi}^{v,k}]_\tau \qquad (\pi, v) \in \mathcal I_\tau^{k+1}.
\end{equation}
In the normal direction, this enforces the condition that the negative and positive
subfaces must be in contact in the next iteration $k+1$. In tangential direction the
enforced condition is dependent on the previous Newton iteration. If the subface sliding
velocity was zero in the previous iteration, $[\dot{\vec u}_\pi^{v, k}]=\vec 0$,
Equation~\eqref{eq:stick} enforces the condition that the tangential velocity is zero in
the next iteration, $[\dot{\vec u}_\pi^{v, k+1}]=\vec 0$. If the subface sliding
velocity was different from zero in the previous iteration, the Newton update does not
immediately enforce zero sliding velocity, however, as the algorithm converges we have
$F_\pi^v\lambda_{\pi n}^{v,k+1} = b_\pi^{v,k}$, and the sliding velocity,
$[\dot{\vec u}_\pi^k]$, for the sticking subfaces is set to zero.  For subfaces in
contact and sliding, we enforce
\begin{equation}\label{eq:sliding}
  \begin{aligned}
    {}[\vec u_\pi ^{v,k+1}]_n &= g_\pi^v, \\
    \vec\lambda_{\pi\tau}^{v,k+1} + L_\pi^{v,k}[\dot{\vec u}_{\pi}^{v,k+1}]_\tau +
    F_\pi^v\vec v_\pi^{v,k}\lambda_{\pi n}^{v,k+1}&= \vec r_\pi^{v,k} + b_\pi^{v,k}\vec
    v_f^k,
  \end{aligned}
  \qquad (\pi, v) \in \mathcal A^{k+1}.
\end{equation}
Again, this enforces the condition that the negative and positive subfaces be in contact
at the next iteration $k+1$. In the tangential direction, the condition approximates the
sliding direction and distance. The matrices and vectors are:
\begin{align}
  L_\pi^{v,k} &= c (( I_{d-1} - M_\pi^{v,k})^{-1} - I_{d-1}) \label{eq:robin_weight}\\
  \vec v_\pi^{v,k} &= ( I_{d-1} - M_\pi^{v,k})^{-1} 
                     \frac{-\vec\lambda_{\pi\tau}^{v,k} 
                     + c\vec
                     [\dot{\vec u}_\pi^{v,k}]_\tau}{\lVert -\vec\lambda_{\pi\tau}^{v,k} + 
                     c\vec [\dot{\vec u}_\pi^{v,k}]_\tau\rVert} \nonumber\\
  \vec r_\pi^{v,k} &= -(I_{d-1} - M_\pi^{v,k})^{-1} e_\pi^{v,k} 
                     Q_\pi^{v,k}(-\vec\lambda_{\pi\tau}^{v,k}
                     + c[\dot{\vec u}_\pi^{v,k}]_\tau),\nonumber
\end{align}
where $I_{d-1}$ is the $(d-1\times d-1)$ identity matrix and
$M_\pi^{v,k} = e_\pi^{v,k}(I_{d-1} - Q_\pi^{v,k})$ with
\begin{align*}
  Q_\pi^{v,k} = \frac{-\vec\lambda_{\pi\tau}^{v,k}(-\vec\lambda_{\pi\tau}^{v,k} + 
  c\vec
  [\dot{\vec u}_\pi^{v,k}]_\tau)^\top}{b_\pi^{v,k}\lVert -\vec\lambda_{\pi\tau}^{v,k} 
  + c\vec
  [\dot{\vec u}_\pi^{v,k}]_\tau\rVert},\qquad
  e_\pi^{v,k} = \frac{b_\pi^{v,k}}{\lVert-\vec
  \lambda_{\pi\tau}^{v,k} + c\vec [\dot{\vec u}_\pi^{v,k}]_\tau\rVert}.
\end{align*}

\subsubsection*{Regularization}
For the subfaces in the inactive set $\mathcal I_n^{k+1}$, i.e., the subfaces not in
contact, the Newton update gives a homogeneous Neumann boundary condition. For the
subfaces in the contact sets $\mathcal I_\tau^{k+1}$ and $\mathcal A^{k+1}$, the Newton
update gives a Dirichlet condition in the normal direction and a Robin boundary
condition in the tangential direction. This Robin condition guarantees positive
definiteness of the system only if $L_\pi^{v,k}$, defined by
Equation~\eqref{eq:robin_weight}, is positive definite. In the converged limit, the
matrix $L_\pi^{v,k}$ is a positive definite
matrix~\cite{hueber2008primal-dual}. However, during the iterations, there is no
guarantee that this will hold. We therefore add a regularization to the Robin conditions
by replacing $Q_\pi^{v,k}$ by
\begin{equation*}
  \tilde{Q}_\pi^{v,k} =  
  \frac{-\vec\lambda_{\pi\tau}^{v,k}(-\vec\lambda_{\pi\tau}^{v,k} + c\vec
    [\dot{\vec u}_\pi^{v,k}]_\tau)^\top}{\max(b_{\pi}^{v,k}, \lVert\vec\lambda_{\pi\tau}^{v,k}\rVert)\lVert
    -\vec\lambda_{\pi\tau}^{v,k} + c\vec
    [\dot{\vec u}_\pi^{v,k}]_\tau\rVert},
\end{equation*}
which is only different from $Q_\pi^{v,k}$ when the inequalities in
Equation~\eqref{eq:discrete_coulomb_friction} are violated. Further, we define
\begin{equation*}
  \alpha_\pi^{v,k} = \frac{(-\vec
    \lambda_{\pi\tau}^{v,k})^\top(-\vec\lambda_{\pi\tau}^{v,k} + c[\dot{\vec u}_\pi^{v,k}]_\tau)}
  {\lVert\vec\lambda_{\pi\tau}^{v,k}\rVert\lVert 
    -\vec\lambda_{\pi\tau}^{v,k} 
    +c\vec
    [\dot{\vec u}_\pi^{v,k}]_\tau\rVert},\qquad \delta_\pi^{v,k} = \min\left(\frac{\lVert
      \vec\lambda_{\pi\tau}^{v,k}\rVert}{\lambda_{\pi n}^{v,k}}, 1 \right),
\end{equation*}
and
\begin{equation*}
  \beta_\pi^{v,k} =
  \begin{cases}
    \frac{1}{1-\alpha_\pi^{v,k}\delta_\pi^{v,k}}, & \text{if } \alpha_\pi^{v,k} <0\\
    1, &\text{otherwise}.
  \end{cases}
\end{equation*}
Using the notation that tilde ($\tilde{\cdot}$) denotes the regularization, we have
$\tilde{M}_\pi^{v,k} = e_\pi^{v,k}(I_{d-1} - \tilde{Q}_\pi^{v,k})$ and replace the
matrix $L_\pi^{v,k}$ from Equation~\eqref{eq:robin_weight} by
\begin{equation*}
  \tilde{L}_\pi^{v,k} = c((I_{d-1} - \beta_\pi^{v,k}\tilde{M}_\pi^{v,k}) - I_{d-1}),
\end{equation*}
which guarantees its positive definiteness~\cite{hueber2008primal-dual}. As the iterates
$(\vec u^k, \vec \lambda^k)$ converge to the solution, the regularization
$\tilde{Q}_\pi^{v,k}\rightarrow Q_\pi^{v,k}$,
$\tilde{L}_\pi^{v,k}\rightarrow L_\pi^{v,k}$, and $\beta_\pi^{v,k}\rightarrow 1$, and we
obtain the original system of equations.

\subsection{Discrete system of equations}
We end Section~\ref{sec:discretization} with a summary of the discrete system of
equations that is solved at each Newton iteration, and we state the discrete linearized
version of Equations~\eqref{eq:biot}-\eqref{eq:coulomb_friction} as:
\begin{equation}
  \label{eq:biot_discrete}
  \begin{aligned}
    A \vec u  + B \vec p + C \vec \lambda &= \vec b_u,\\
    D\vec u + E\vec p + F\vec \lambda &= \vec b_p,\\
    G\vec u + H \vec p + J\vec \lambda &= \vec r.
  \end{aligned}
\end{equation}
The first line is the discrete momentum balance, and the matrices $A$, $B$, $C$, and the
vector $\vec b_u$ are obtained by considering Equation~\eqref{eq:momentum_balance_cell}
for all cells and assembling the coefficients in the global matrices. Similarly, the
second row corresponds to the discrete flux balance, and the equation is obtained by
considering Equation~\eqref{eq:flux_balance_cell} over all cells and assembling the
coefficients in the global matrices. The matrix $F$ appears due to the dependence on
$\vec \lambda$ in the local systems for the displacement gradients, as given by
Equation~\eqref{eq:gradu}. It is worth pointing out that the matrices $A$, $B$, $D$, and
$E$ are the same matrices as are obtained by the finite-volume scheme in a poroelastic
domain without any fractures~\cite{nordbotten2016stable}. The last row of
Equation~\eqref{eq:biot_discrete} corresponds to the linearization of the complementary
functions~\eqref{eq:comp_n} and~\eqref{eq:comp_t}, and the matrices are obtained by
assembling Equations~\eqref{eq:open}-\eqref{eq:sliding} for each subface on the
fracture. The dependence of the pressure in the contact law, given by the matrix $H$, is
due to the pressure dependence on the poroelastic stress (see second row of
Equation~\ref{eq:biot}) as well as the pressure dependence on the displacement gradients
given by Equation~\eqref{eq:gradu}.

From a computational point of view, it is worth noting that during the Newton iteration,
only the matrices $G$, $H$, $J$, and the vector $\vec r$ will change. This means that
updating the discretization is inexpensive as it is only a local update for the subcells
bordering the fractures.

\section{Numerical examples}
Four numerical examples are given. For the first two, we neglect the fluid contribution
to the mechanical stress to investigate the performance of the numerical approach for
the purely mechanical contact problem, i.e., we set $\alpha=0$. In all of the examples,
Young's modulus is $E_0 = 4$~GPa, the Poisson ratio is $\nu = 0.2$, and the initial gap
of the fractures is $g=0$. In our experience, the algorithm is quite robust with respect
to the numerical parameter $c$, and in the examples, it is fixed to $c=100$~GPa/m.

We assign a space varying coefficient of friction so that the slip of the fractures will
arrest before it reaches the fracture tips. This choice of the friction coefficient is
done to obtain a solution with high enough regularity to study the convergence in
stress. If the slip of the fractures reaches the fracture tips, the solution will
contain singularities in the stress, which reduces the regularity of our solution. Note
that our method is not restricted to the regularized solution, as discussed more
thoroughly in Appendix~\ref{sec:appendix}.

The discrete solution is denoted $\vec u_h$, which is interpreted as the piecewise
constant function over each cell $K\ni\vec x$ such that $\vec u_h(\vec x) = \vec
u_K$. The discrete solution $\vec\lambda_h$ for the Lagrange multiplier is defined as
piecewise constant on each face $\pi$ on $\Gamma^+$ and is equal to the area weighted
sum of the subface values,
$\vec\lambda_h(\vec x)m_\pi = \sum_{v\in\mathcal V_\pi}m_\pi^v \lambda_\pi^v$,
$\vec x\in\pi$. Likewise, the displacement jump is defined as the piecewise constant on
each face, $\pi$, on $\Gamma^+$ corresponding to the subface average,
$[\vec u_h(\vec x)] = \frac{1}{|\mathcal V_\pi|}\sum_{v\in\mathcal V_\pi}[\vec
u_\pi^v]$, $\vec x\in\pi$, where $|\mathcal V_\pi|$ is the number of subfaces of the
face, $\pi$, which is equal to three if $\pi$ is a triangle. The continuous solution is
denoted by the pair $(\vec u, \vec\lambda)$.

We define the relative error of a discrete variable $\xi_h$ in a domain $\gamma$ as
\begin{equation}
  \label{eq:relative_error}
  \varepsilon_\gamma(\xi_h, \xi) = \frac{\lVert \vec \xi_h - \vec 
    \xi\lVert_{\gamma}}{\lVert \vec
    \xi\rVert_{\gamma}},
\end{equation}
where $\xi$ is a reference solution and $||\cdot||_\gamma$ is the $L_2$ norm over the
domain $\gamma$. The Newton iteration is terminated when the change in the solution is
below a given stopping criterion:
\begin{equation}\label{eq:newton_error}
  \varepsilon_\Omega(\vec u_h^{k+1}, \vec u_h^{k}) < \delta,
\end{equation}
where $k$ is the Newton iteration index.

To solve the linear system of equations at each Newton iteration, a direct solver is
used if the number of degrees of freedom is less than 10 000, else, an iterative solver
is used. The iterative solver uses a preconditioned GMRES iteration that is based on a
Schur complement strategy, where the mechanics-fluid subsystem is approximated by a
single fixed stress iteration, see~\cite{kim2011,both2017} for details. Within the fixed
stress iteration, the mechanics problem is solved by one AMG iteration, as implemented
in~\cite{oson2018}, while a direct solver is applied to the flow problem. In the simpler
case of a pure mechanics problem, the fixed stress iteration is replaced with an AMG
iteration on the mechanics subproblem.

\begin{figure}[h]
  \begin{subfigure}{0.49\textwidth}
    \centering \def\svgwidth{\textwidth}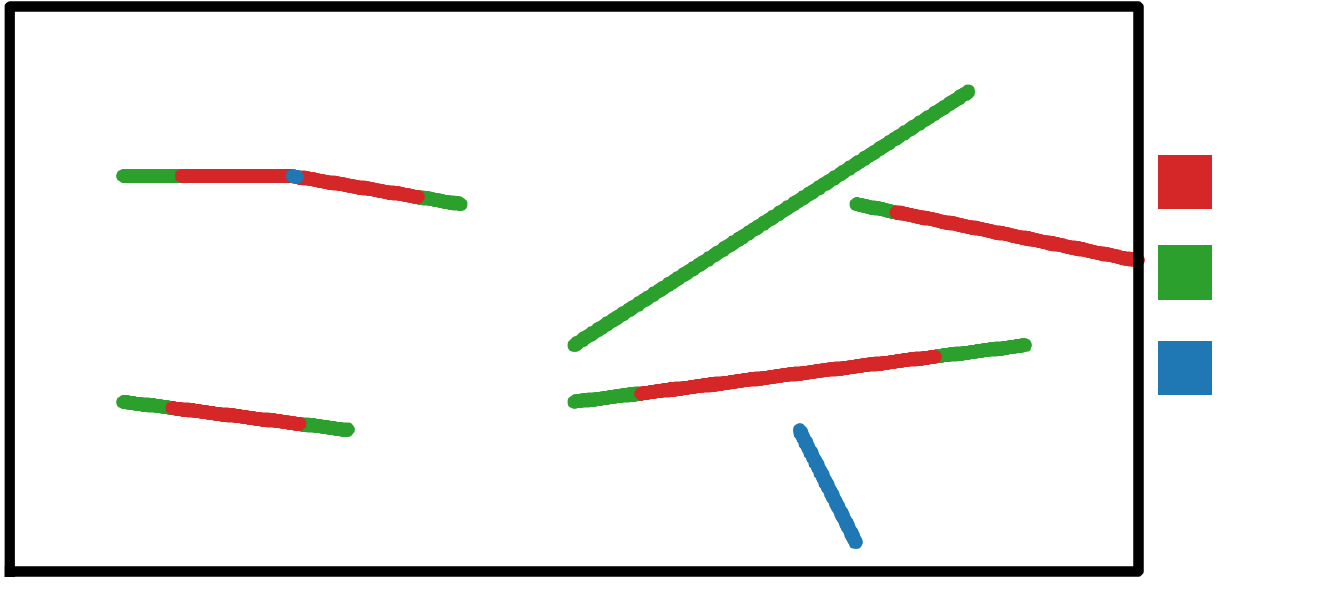
  \end{subfigure}
  \hfill
  \begin{subfigure}{0.43\textwidth}
    \fontsize{8pt}{10pt}\selectfont
    \def\svgwidth{\textwidth}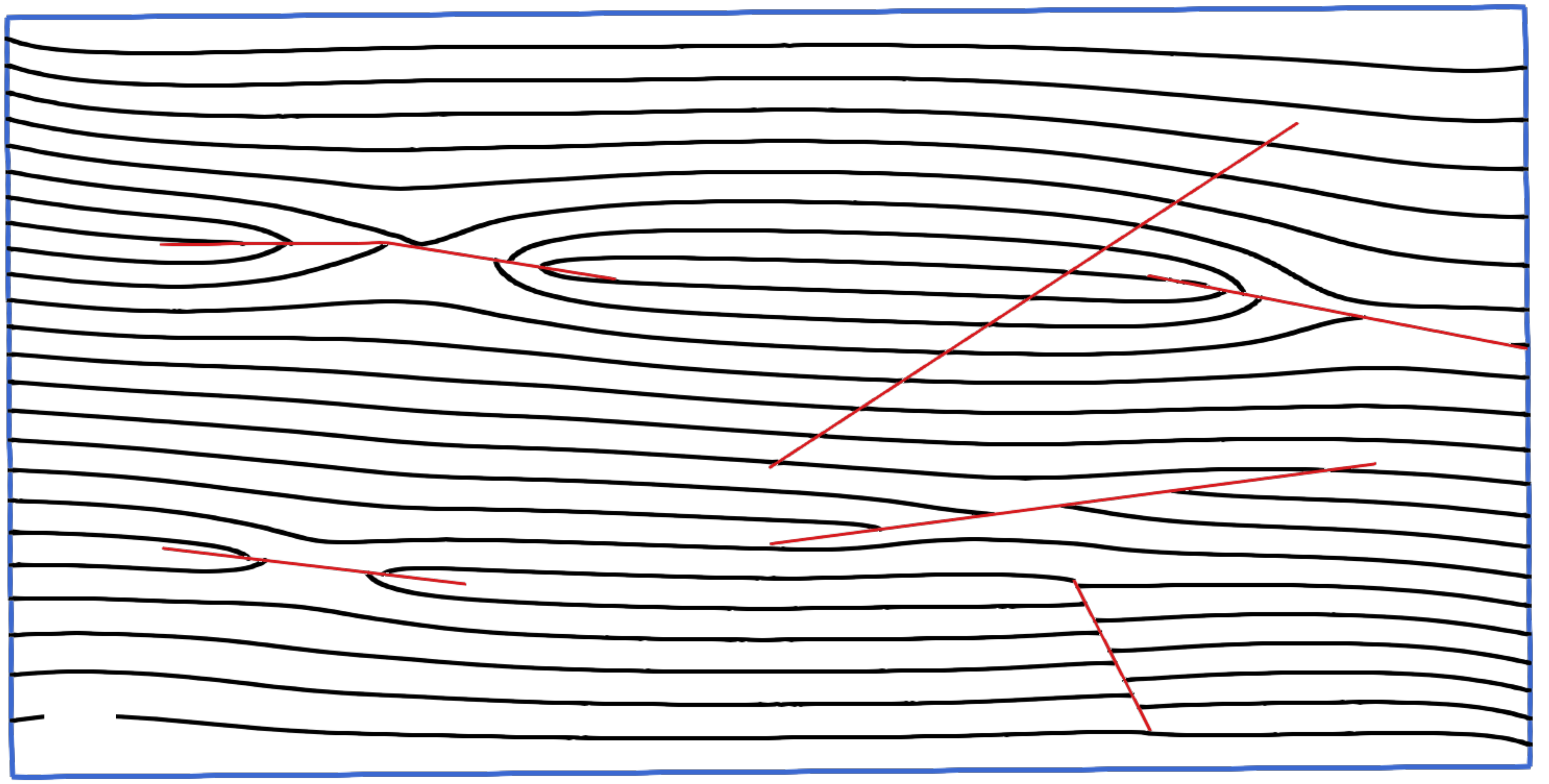
  \end{subfigure}
  \caption{The fractured domain from Example 1. Left: The black box corresponds to the
    domain boundaries, while the fractures, labeled by a number, are represented by
    lines.  The fractures are colored in three colors; segments that slip are red,
    segments that stick are green, and segments that open are blue. Right: Contour plot
    of the $x$-component of the displacement $\vec u$.  The labels on the contours are
    given in millimeters. The red lines represent the fractures.\label{fig:slip_stick}}
\end{figure}

The computer code has been implemented in the open source Python toolbox
PorePy~\cite{keilegavlen20017porepy}, which has an interface for meshing in
Gmsh~\cite{geuzaine2009gmsh}. The run scripts for the examples are open
source~\cite{berge2019finite}. ParaView~\cite{ayachit2015paraview} was used to make
Figures~\ref{fig:slip_stick} and~\ref{fig:ex2_mortar_grid_lam_u}.

\vspace{-2em}
\subsection{Example 1}
The first example is a domain $2$~m~$\times\ 1$~m with six fractures, as depicted in
Figure~\ref{fig:slip_stick}. This example includes difficult cases such as a fracture
with a kink and a fracture reaching the boundary. An advantage of our finite volume
method is that no special treatment is needed to handle these cases because the degrees
of freedom are located in the cell and subface centers and not on the nodes. In this
example, we do not consider any fluid and solve only for the linear elasticity in
Equation~\eqref{eq:biot} coupled to the contact conditions given in
Equations~\eqref{eq:non_penetration} and~\eqref{eq:coulomb_friction}. The bottom
boundary is fixed, the two vertical boundaries are free, and at the top boundary a
Dirichlet condition $\vec g_{u,D} = [0.005, -0.002]^\top$~m is assigned. The initial
guess in the Newton iteration is $\vec u=\vec 0$~m, $\lambda_n =-100$~Pa and
$\vec\lambda_\tau=\vec 0$~Pa, i.e., zero displacement and all fractures in contact and
sticking. The coefficient of friction is for each fracture $i=1\ldots6$ set to
$F_i(\vec x) = 0.5(1 + \exp(-D_i(\vec x)^2 / 0.005\text{m}^2),\ \vec x\in \Gamma_i^+$,
where $D_i(\vec x)$ is the distance from $\vec x$ to the tips of fracture $i$. Note that
the bend in Fracture $1$ and the right end of Fracture $5$ are not considered tips for
the distance function $D$, and thus the coefficient of friction at these points is
$F\approx 0.5$.

\begin{figure}
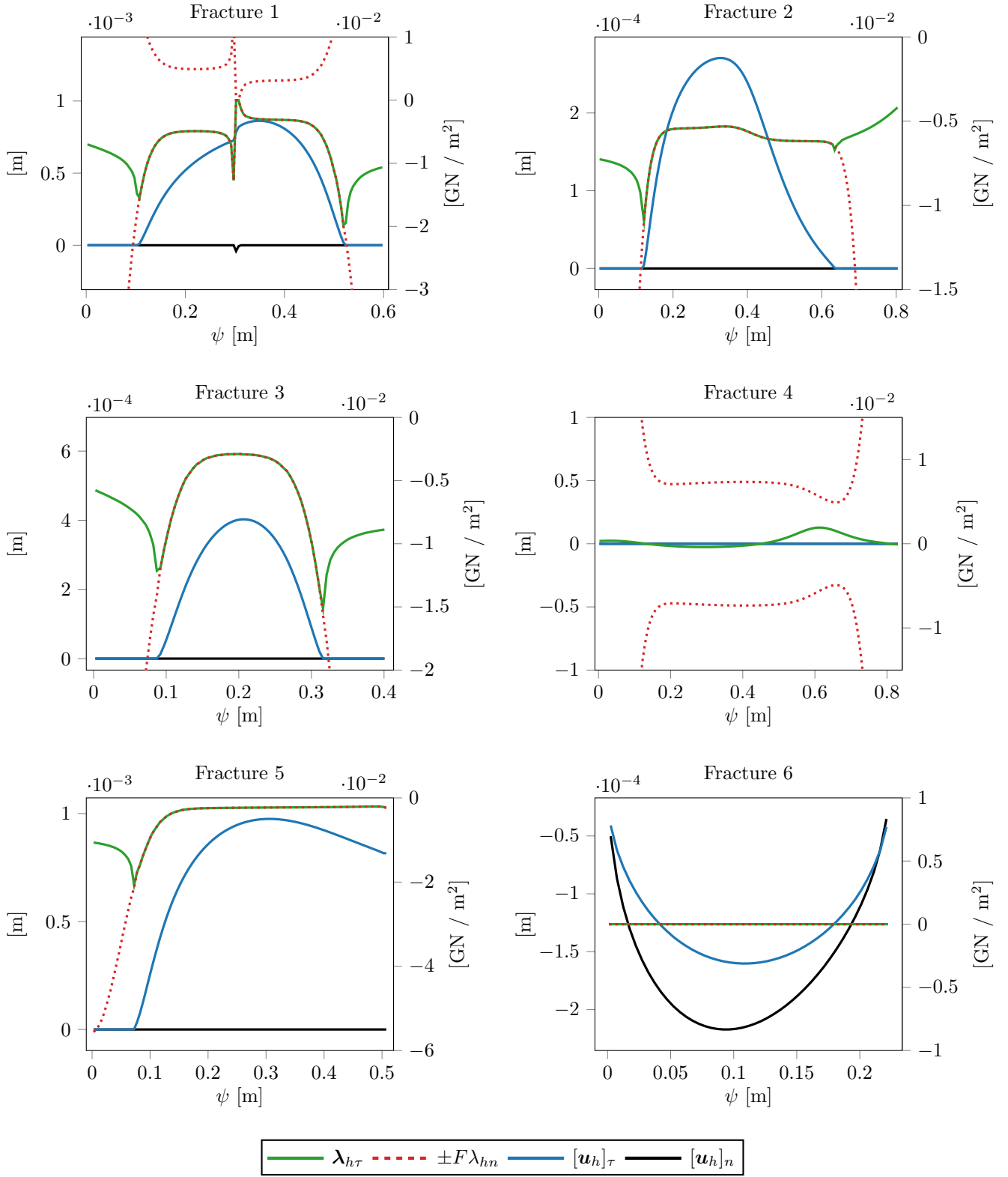

  \centering \def\figwidth{0.4\textwidth} \def\hfig{0.1\textwidth}
  \hspace*{-5em}\begin{subfigure}{\figwidth}
    \def\svgwidth{\textwidth} \import{fig/example1/}{traction_reg0.tex}
  \end{subfigure}
  \hspace{\hfig}
  \begin{subfigure}{\figwidth}
    \def\svgwidth{\textwidth} \import{fig/example1/}{traction_reg1.tex}
  \end{subfigure}
  \hspace*{-5em}\begin{subfigure}{\figwidth}
    \def\svgwidth{\textwidth} \import{fig/example1/}{traction_reg2.tex}
  \end{subfigure}
  \hspace{\hfig}
  \begin{subfigure}{\figwidth}
    \def\svgwidth{\textwidth} \import{fig/example1/}{traction_reg3.tex}
  \end{subfigure}
  \hspace*{-5em}\begin{subfigure}{\figwidth}
    \def\svgwidth{\textwidth} \import{fig/example1/}{traction_reg4.tex}
  \end{subfigure}
  \hspace{\hfig}
  \begin{subfigure}{\figwidth}
    \def\svgwidth{\textwidth} \import{fig/example1/}{traction_reg5.tex}
  \end{subfigure}
  \begin{subfigure}{\textwidth}
    \centering
    \begin{tikzpicture}
      \definecolor{blue}{rgb}{0.12156862745098,0.466666666666667,0.705882352941177}
      \definecolor{green}{rgb}{0.172549019607843,0.627450980392157,0.172549019607843}
      \definecolor{red}{rgb}{0.83921568627451,0.152941176470588,0.156862745098039}
      \path[ultra thick, draw, color=green] (0,0)--(1,0) node[right, black] (A)
      {$\vec\lambda_{h\tau}$}; \path[ultra thick, draw, dashed, color=red]
      (A.east)--($(A.east)+(1,0)$) node[right, black] (B) {$\pm F\lambda_{hn}$};
      \path[ultra thick, draw, color=blue] (B.east)--($(B.east) + (1,0)$) node[right,
      black] (C) {$[\vec u_{h}]_\tau$}; \path[ultra thick, draw]
      (C.east)--($(C.east) + (1,0)$) node[right] (D) {$[\vec u_{h}]_n$}; \draw[thick]
      ($(A.south west) + (-1.1, 0)$) rectangle (D.north east);
    \end{tikzpicture}
  \end{subfigure}
  \caption{Results from Example 1 showing $\vec\lambda_h$ and $[\vec u_h]$ for the
    reference solution that has an average of 103 faces along each fracture. The shear
    component of the Lagrange multiplier $\vec\lambda_{h\tau}$ (green), friction bound
    $\pm F\lambda_{hn}$ (red dashes), tangential displacement jump $[\vec u_h]_\tau$
    (blue), and normal displacement jump $[\vec u_h]_n$ (black) for all fractures. The
    $x$-axis shows the distance $\psi$ from the leftmost end of the fracture. The
    subplots are arranged from top left to bottom right according to the fracture number
    given in Figure~\ref{fig:slip_stick}.  \label{fig:ex1_traction}}
\end{figure}

A contour plot of the solution is shown in Figure~\ref{fig:slip_stick} where the
discontinuous displacement over the sliding or opening fractures can clearly be seen. To
better visualize the different behaviors of the fractures, the fracture regions that are
slipping, sticking, and opening are plotted in different colors in
Figure~\ref{fig:slip_stick}. For Fracture 1, the top boundary is sliding to the right,
while the bottom boundary is sliding to the left. This situation causes the fracture to
open in a small segment after the bend.  Figure~\ref{fig:ex1_traction} shows the shear
component of the Lagrange multiplier as well as the friction bound and displacement
jump. At the bend of fracture 1, there is a singularity in the stress that causes the
sharp increase in the Lagrange multiplier. For Fracture 2, we observe a change in the
shear and normal component of the Lagrange multiplier at approximately the midpoint that
is caused by the opening of Fracture 6. In the vicinity of the fracture tips, there is a
sharp increase in the shear component of the Lagrange multiplier as the fractures change
behavior from sliding to sticking.

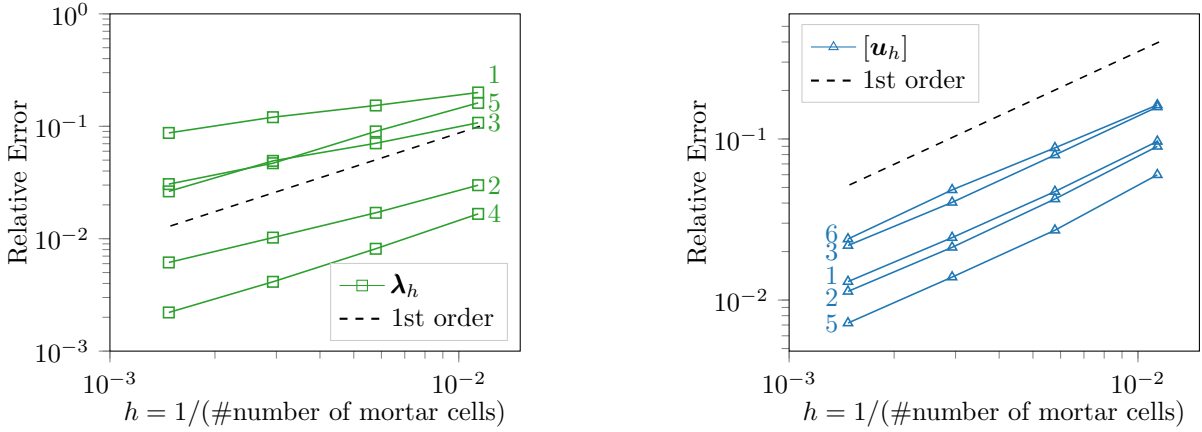
\begin{figure}
  \def\figwidth{0.4\textwidth}
  \begin{subfigure}[t]{\figwidth}
    \centering \def\svgwidth{\textwidth}
    \begin{tikzpicture}

\definecolor{green}{rgb}{0.172549019607843,0.627450980392157,0.172549019607843}

\begin{axis}[
width=\svgwidth,
legend cell align={left},
legend entries={{$\pmb \lambda_h$},{1st order}},
legend style={at={(0.97,0.03)}, anchor=south east, draw=white!80.0!black},
tick align=outside,
tick pos=left,
x grid style={white!69.01960784313725!black},
xlabel={$h = 1 / (\#${number of mortar cells}$)$},
xmin=1e-3, xmax=1.5e-2,
ymin=1e-3, ymax=1e0,
xmode=log,
y grid style={white!69.01960784313725!black},
ylabel={Relative Error},
ymode=log
]
\addlegendimage{mark=square, green}
\addlegendimage{thick, black, dashed}
\addplot [semithick,dashed]
table [row sep=\\]{%
0.0114942528735632	0.1 \\
0.00148367952522255 	0.012908011869436208 \\
};
\addplot [semithick, green, mark=square, mark size=2, mark options={solid,fill=black}]
table [row sep=\\]{%
0.0113636363636364	0.19976067216498 \\
0.00578034682080925	0.153163078811876 \\
0.00293255131964809	0.120458327729829 \\
0.00147492625368732	0.0873138566016368 \\
} node [above right,pos=0] {$1$};
\addplot [semithick, green, mark=square, mark size=2, mark options={solid,fill=black}, forget plot]
table [row sep=\\]{%
0.0113636363636364	0.0299184637515252 \\
0.00578034682080925	0.0169989124622728 \\
0.00293255131964809	0.010227778286855 \\
0.00147492625368732	0.00614880275917388 \\
} node [right,pos=0] {$2$};
\addplot [semithick, green, mark=square, mark size=2, mark options={solid,fill=black}, forget plot]
table [row sep=\\]{%
0.0113636363636364	0.107507794487453 \\
0.00578034682080925	0.0705759836387564 \\
0.00293255131964809	0.0492525114110022 \\
0.00147492625368732	0.0263692066055125 \\
} node [right,pos=0] {$3$};
\addplot [semithick, green, mark=square, mark size=2, mark options={solid,fill=black}, forget plot]
table [row sep=\\]{%
0.0113636363636364	0.0166305356316782 \\
0.00578034682080925	0.00813966349100972 \\
0.00293255131964809	0.00413686782908446 \\
0.00147492625368732	0.00220718179735522 \\
} node [right,pos=0] {$4$};
\addplot [semithick, green, mark=square, mark size=2, mark options={solid,fill=black}, forget plot]
table [row sep=\\]{%
0.0113636363636364	0.161567416231477 \\
0.00578034682080925	0.0900348427517522 \\
0.00293255131964809	0.0469404172700777 \\
0.00147492625368732	0.0305936465074972 \\
} node [right,pos=0] {$5$};
\end{axis}

\end{tikzpicture}
  \end{subfigure}
  \hspace{0.1\textwidth}
  \begin{subfigure}[t]{\figwidth}
    \centering \def\svgwidth{\textwidth}
    \begin{tikzpicture}
\definecolor{blue}{rgb}{0.12156862745098,0.466666666666667,0.705882352941177}
\begin{axis}[
width=\svgwidth,
legend cell align={left},
legend entries={{$[\vec u_h]$},{1st order}},
legend style={at={(0.03,0.97)}, anchor=north west, draw=white!80.0!black},
tick align=outside,
tick pos=left,
x grid style={white!69.01960784313725!black},
xlabel={$h = 1 / (\#${number of mortar cells}$)$},
xmin=1e-3, xmax=1.5e-2,
xmode=log,
y grid style={white!69.01960784313725!black},
ylabel={Relative Error},
ymode=log
]
\addlegendimage{blue, mark=triangle, mark size = 2}
\addlegendimage{thick, black, dashed}
\addplot [semithick,dashed]
table [row sep=\\]{%
0.0114942528735632	0.4 \\
0.00148367952522255 	0.05163204747774483 \\
};
\addplot [semithick, blue, mark=triangle, mark size=2]
table [row sep=\\]{%
0.0113636363636364	0.096610945517093 \\
0.00578034682080925	0.0471497899576949 \\
0.00293255131964809	0.0244147636706353 \\
0.00147492625368732	0.0130174339611638 \\
} node [above=0.5ex, left, pos=1] {$1$};
\addplot [semithick, blue, mark=triangle, mark size=2]
table [row sep=\\]{%
0.0113636363636364	0.0899647978914078 \\
0.00578034682080925	0.0424120176431605 \\
0.00293255131964809	0.0212473866737652 \\
0.00147492625368732	0.0113166946847267 \\
} node [below=0.5ex, left, pos=1] {$2$};
\addplot [semithick, blue, mark=triangle, mark size=2]
table [row sep=\\]{%
0.0113636363636364	0.157970547924419 \\
0.00578034682080925	0.0795742267604959 \\
0.00293255131964809	0.0403979034006932 \\
0.00147492625368732	0.0218108638753825 \\
} node [below=0.5ex, left,pos=1] {$3$};

\addplot [semithick, blue, mark=triangle, mark size=2]
table [row sep=\\]{%
0.0113636363636364	0.0599170411794384 \\
0.00578034682080925	0.0271341927913356 \\
0.00293255131964809	0.0138878442228929 \\
0.00147492625368732	0.0072031062712245 \\
} node [left,pos=1] {$5$};
\addplot [semithick, blue, mark=triangle, mark size=2]
table [row sep=\\]{%
0.0113636363636364	0.162559939005485 \\
0.00578034682080925	0.0880535616593204 \\
0.00293255131964809	0.0483185864672939 \\
0.00147492625368732	0.0239351228475929 \\
} node [above=0.5ex, left,pos=1] {$6$};
\end{axis}

\end{tikzpicture}
  \end{subfigure}
  \caption{Convergence rates for the Lagrange multiplier $\vec\lambda_h$ (left) and the
    displacement jump $[\vec u_h]$ (right) for each separate fracture in Example 1.  The
    error is measured as the relative errors
    $\varepsilon_{\Gamma_i^+}([\vec u_h], [\vec u])$ and
    $\varepsilon_{\Gamma_i^+}(\vec\lambda_h, \vec\lambda)$ for each fracture
    $\Gamma_i^+$. The line numbering corresponds to the fracture numbers given in
    Figure~\ref{fig:slip_stick}. \label{fig:ex1_rate_mortar}}
\end{figure}

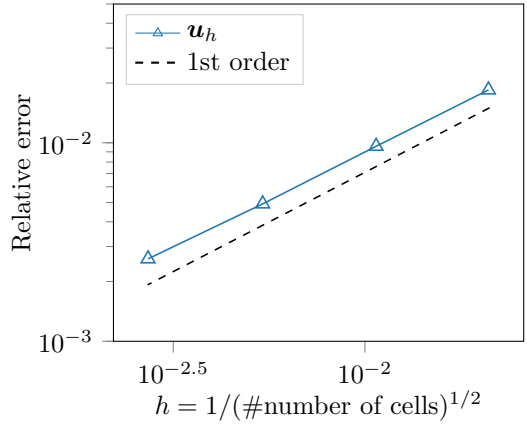
\begin{figure}
  \centering \def\svgwidth{0.4\textwidth}
  \begin{tikzpicture}
\definecolor{blue}{rgb}{0.12156862745098,0.466666666666667,0.705882352941177}
  \begin{axis}[
    width=\svgwidth,
legend cell align={left},
legend entries={{$\vec u_h$},{1st order}},
legend style={at={(0.03,0.97)}, anchor=north west, draw=white!80.0!black},
tick align=outside,
tick pos=left,
x grid style={white!69.01960784313725!black},
xlabel={$h = 1 / (\#${number of  cells}$)^{1/2}$},
xmode=log,
y grid style={white!69.01960784313725!black},
ylabel={Relative error},
ymode=log,
ymin=0.001,
ymax=0.05
]
\addlegendimage{mark=triangle, blue}
\addlegendimage{thick, black, dashed}
\addplot [semithick, dashed]
table [row sep=\\]{%
0.0211335742104886	0.015 \\
0.00271083353596815	0.0019240712732511397 \\
};
\addplot [semithick, blue, mark=triangle, mark size=3, mark options={solid,fill=black}]
table [row sep=\\]{%
0.0209933986140406	0.0184589647683222 \\
0.0106959518669558	0.00960149529139043 \\
0.00541014976448614	0.00492223508550284 \\
0.00271785293542299	0.0026026093151755 \\
};
\end{axis}

\end{tikzpicture}
  \caption{Convergence rates for the cell-centered displacement in Example 1. The error
    is the relative error $\varepsilon_{\Omega}(\vec u_h, \vec u)$, as defined by
    Equation~\eqref{eq:relative_error}.\label{fig:ex1_rate}}
\end{figure}

As a reference solution $(\vec u,\vec\lambda)$, we use the solution calculated for a
fine mesh using 1.7 million degrees of freedom. The second finest mesh has 270 thousand
degrees of freedom and is the mesh used for the results in Figure~\ref{fig:ex1_traction}
and~\ref{fig:slip_stick} . In Figure~\ref{fig:ex1_rate_mortar}, the relative errors
$\varepsilon_{\Gamma_i^+}([\vec u_h], [\vec u])$ and
$\varepsilon_{\Gamma_i^+}(\vec\lambda_h, \vec\lambda)$, given by
Equation~\eqref{eq:relative_error}, are plotted for each fracture $i=1,\ldots,6$. For
the displacement jump, the convergence is of first-order for all fractures except
Fracture 4, which is correctly predicted to be sticking (and thus, the error is
zero). For the Lagrange multiplier $\vec\lambda_h$, we observe first-order convergence
for Fractures $4$ and $5$, while the error for Fracture $6$ is zero. The convergence
rates for traction is typically observed to be of first-order in the $L_2$ norm and
second-order in the $2$-norm for the finite-volume
scheme~\cite{keilegavlen2017finite,nordbotten2016stable}, but the Lagrange multiplier
for fractures $1$, $2$ and $3$ shows somewhat lower convergence rates than
first-order. However, this is not surprising due to the low regularity of the Lagrange
multipliers. Figure~\ref{fig:ex1_rate} shows first-order convergence of the error for
the discrete displacement $\vec u_h$ in the 2d domain $\Omega$. Finally,
Table~\ref{table:newton_count} shows that the number of Newton iterations do not grow
significantly when the mesh is refined.

\subsection{Example 2}\label{sec:example2}
In this example $\Omega$ is a 3d domain
$(-200, 300)\times(-200, 300)\times(-300, 300)$~m with two circular fractures
approximated by polygons with 10 vertices. The location and geometry of the fractures
are given in Table~\ref{tab:frac_geom}.
\begin{table}
  \begin{center}
    \caption{Fracture geometry in Example 2 and 3. The strike angle is the rotation from
      x-axis in the x-y-plane defining the strike line. The dip angle is rotation around
      the strike line.\label{tab:frac_geom}}
    \begin{tabular}{lcc}
      \toprule
      & Fracture 1 & Fracture 2 \\ \midrule
      Center & $-[10, 30, 80]^\top$ m & $[15, 60, 80]^\top$ m \\
      Radius & 150 m &  150 m \\
      Strike angle & $81.8^\circ$ & $78.3^\circ$ \\
      Dip angle & $43.9^\circ$ & $47.1^\circ$ \\
      \bottomrule
    \end{tabular}
  \end{center}
\end{table}
As in the previous example, no fluid is included. The bottom boundary is fixed, the four
vertical boundaries are rolling, and at the top boundary, a load is applied downwards by
enforcing a Neumann condition $\vec g_{u,N} = [0, 0, -4.5]^\top$~MPa. The coefficient of
friction is for the two fractures, $i=1,2$:
\begin{equation*}
  F_i(\vec x) = 0.5 \exp\left(\frac{10 \text{m}}{R_i - D_i(\vec x)} 
    -\frac{10 \text{m}}{R_i}\right),
\end{equation*}
where $R_i$ is the radius of fracture $i$ and $D_i(\vec x)$ the distance from the center
of the fracture to $\vec x$.

\begin{figure}
  \centering
  \begin{subfigure}[t]{0.25\textwidth}
    \includegraphics[width=\textwidth]{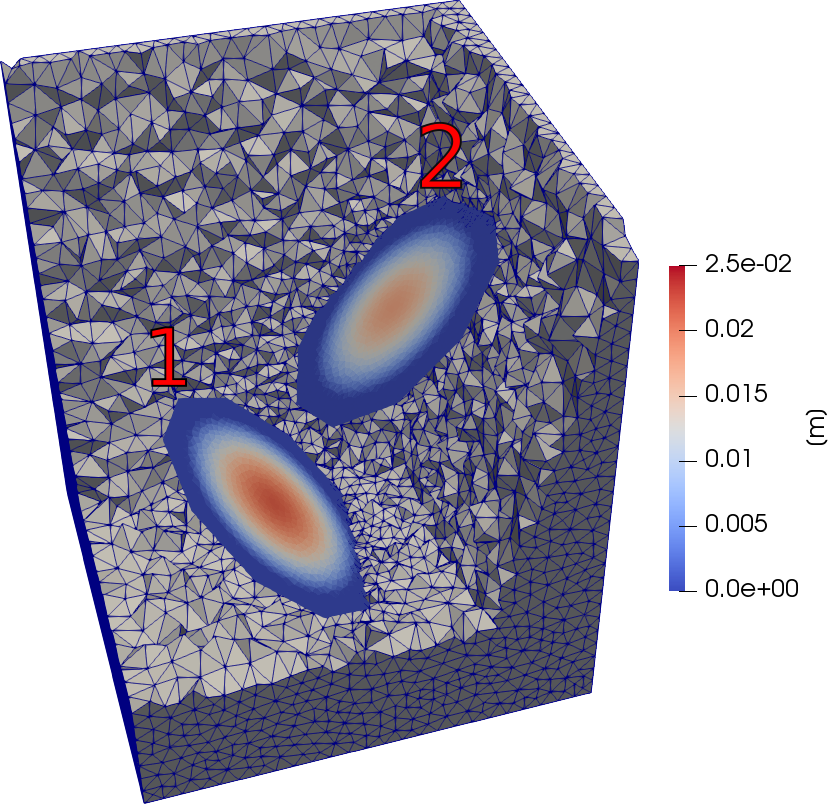}
  \end{subfigure}
  \begin{subfigure}[t]{0.74\textwidth}
    \centering \def\svgwidth{\textwidth}
\begingroup%
  \makeatletter%
  \providecommand\color[2][]{%
    \errmessage{(Inkscape) Color is used for the text in Inkscape, but the package 'color.sty' is not loaded}%
    \renewcommand\color[2][]{}%
  }%
  \providecommand\transparent[1]{%
    \errmessage{(Inkscape) Transparency is used (non-zero) for the text in Inkscape, but the package 'transparent.sty' is not loaded}%
    \renewcommand\transparent[1]{}%
  }%
  \providecommand\rotatebox[2]{#2}%
  \newcommand*\fsize{\dimexpr\f@size pt\relax}%
  \newcommand*\lineheight[1]{\fontsize{\fsize}{#1\fsize}\selectfont}%
  \ifx\svgwidth\undefined%
    \setlength{\unitlength}{799.71661377bp}%
    \ifx\svgscale\undefined%
      \relax%
    \else%
      \setlength{\unitlength}{\unitlength * \real{\svgscale}}%
    \fi%
  \else%
    \setlength{\unitlength}{\svgwidth}%
  \fi%
  \global\let\svgwidth\undefined%
  \global\let\svgscale\undefined%
  \makeatother%
  \begin{picture}(1,0.44938084)%
    \lineheight{1}%
    \setlength\tabcolsep{0pt}%
    \put(0,0){\includegraphics[width=\unitlength,page=1]{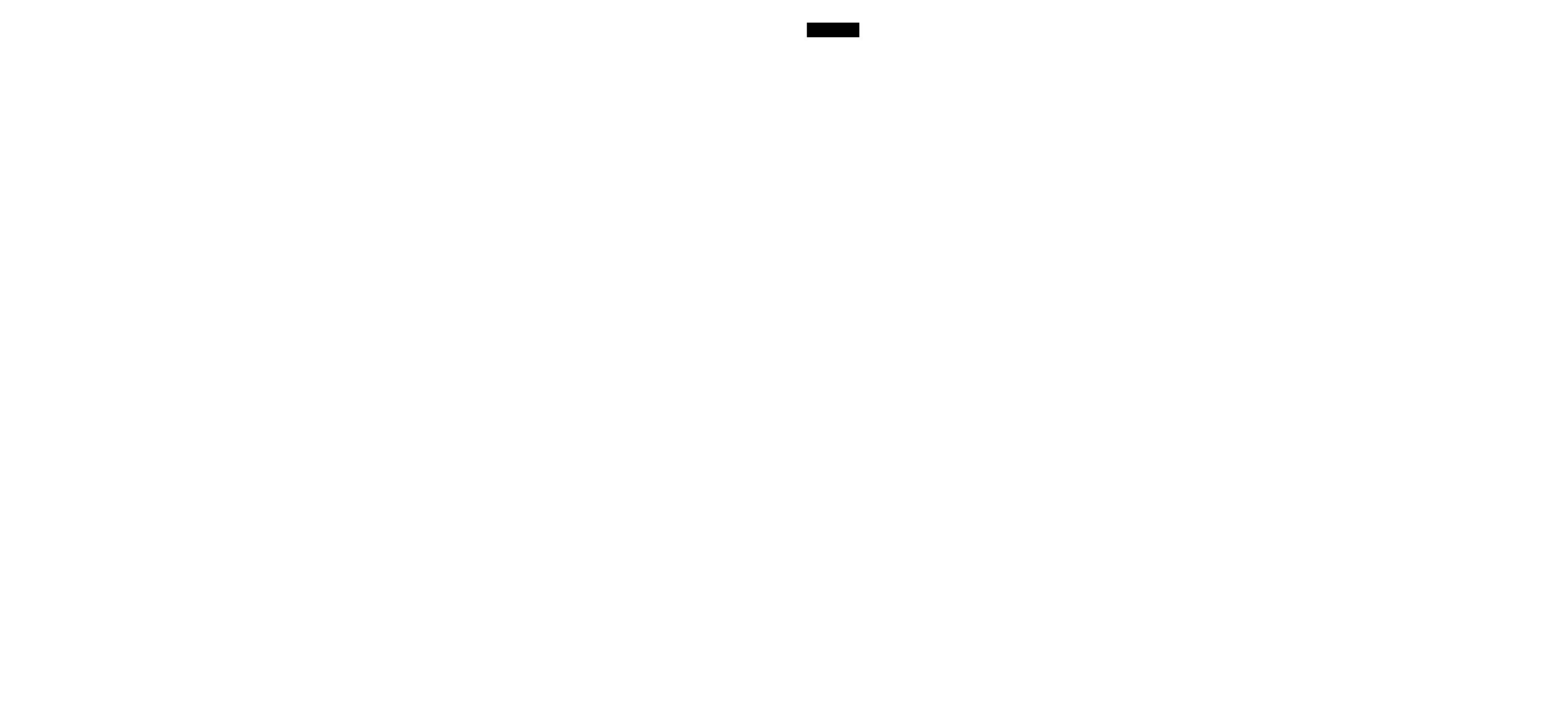}}%
    \put(0.55477103,0.42342208){\color[rgb]{0,0,0}\makebox(0,0)[lt]{\lineheight{1.25}\smash{\begin{tabular}[t]{l}$[\vec u_h]_\tau$\end{tabular}}}}%
    \put(0,0){\includegraphics[width=\unitlength,page=2]{mortar_grid_lam_u_zoom.pdf}}%
    \put(0.55477103,0.39551387){\color[rgb]{0,0,0}\makebox(0,0)[lt]{\lineheight{1.25}\smash{\begin{tabular}[t]{l}$\pmb \lambda_{h\tau}$\end{tabular}}}}%
    \put(0,0){\includegraphics[width=\unitlength,page=3]{mortar_grid_lam_u_zoom.pdf}}%
    \put(0.67610868,0.42427492){\color[rgb]{0,0,0}\makebox(0,0)[lt]{\lineheight{1.25}\smash{\begin{tabular}[t]{l}Contact and sticking\end{tabular}}}}%
    \put(0,0){\includegraphics[width=\unitlength,page=4]{mortar_grid_lam_u_zoom.pdf}}%
    \put(0.67610868,0.39513334){\color[rgb]{0,0,0}\makebox(0,0)[lt]{\lineheight{1.25}\smash{\begin{tabular}[t]{l}Contact and sliding\end{tabular}}}}%
    \put(0,0){\includegraphics[width=\unitlength,page=5]{mortar_grid_lam_u_zoom.pdf}}%
  \end{picture}%
\endgroup%

  \end{subfigure}
  \caption{Results from Example 2 using approximately 250 thousand degrees of
    freedom. Left: Displacement jump $\lVert[\vec u_h]\rVert$ for the two fractures
    indicated by the fracture number. Parts of the 3d mesh are cropped to reveal the
    fractures. Middle: The surface mesh of $\Gamma^+$ of Fracture 1. The red lines show
    the tangential part of the Lagrange multiplier $\vec\lambda_{h\tau}$ while the black
    lines show the scaled displacement jumps $400[\vec u_h]_\tau$. Right: Zoomed view of
    middle figure. \label{fig:ex2_mortar_grid_lam_u}}
\end{figure}

Figure~\ref{fig:ex2_mortar_grid_lam_u} shows the displacement jump $[\vec u_h]_\tau$ and
the shear component of the Lagrange multiplier $\vec\lambda_{h\tau}$. The fractures are
in contact, i.e., the normal displacement jump $[\vec u_h]_n=0$ is zero. Going from two
dimensions to three adds an additional challenge to the contact problem as we have to
find not only the magnitude of the slip but also the direction. The advantage of the
hybrid formulation in combination with a semismooth Newton scheme is that the same
computer code can be used for any dimension, and as observed in the figure, the correct
sliding direction (parallel to the Lagrange multiplier) is found by the algorithm.

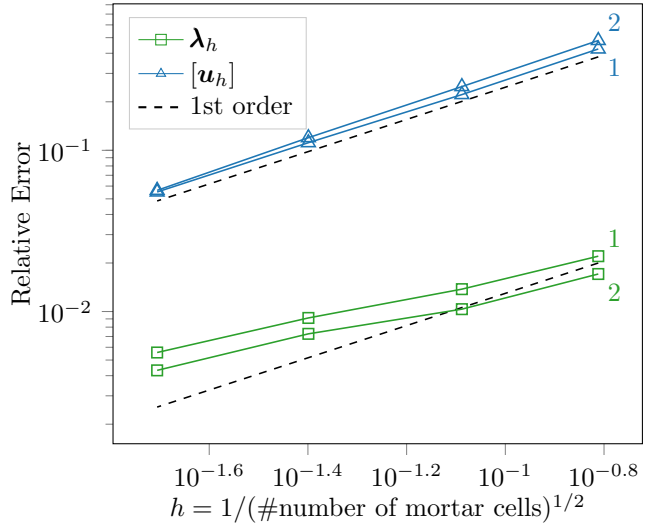
\begin{figure}
  \centering \def\svgwidth{0.49\textwidth}
  \begin{tikzpicture}
  \definecolor{blue}{rgb}{0.12156862745098,0.466666666666667,0.705882352941177}
  \definecolor{green}{rgb}{0.172549019607843,0.627450980392157,0.172549019607843}
\begin{axis}[
    width=\svgwidth,
legend cell align={left},
legend entries={{$\pmb\lambda_h$},{$[\vec u_h]$},{1st order}},
legend style={at={(0.03,0.97)}, anchor=north west, draw=white!80.0!black},
tick align=outside,
tick pos=left,
x grid style={white!69.01960784313725!black},
xlabel={$h = 1 / (\#${number of mortar cells}$)^{1/2}$},
xmode=log,
y grid style={white!69.01960784313725!black},
ylabel={Relative Error},
ymode=log
]
\addlegendimage{mark=square, green}
\addlegendimage{mark=triangle, blue}
\addlegendimage{thick, black, dashed}
\addplot [semithick, dashed]
table [row sep=\\]{%
0.154303349962092	0.02 \\
0.0197027601559775	0.00255386959 \\
};
\addplot [semithick, green, mark=square, mark size=2, mark options={solid,fill=black}]
table [row sep=\\]{%
0.154303349962092	0.0220429913775079 \\
0.0816496580927726	0.0137654339228932 \\
0.0399043442233811	0.00912395795096759 \\
0.0197027601559775	0.00556879838144404 \\
} node [above right, pos=0] {$1$};
\addplot [semithick, green, mark=square, mark size=2, mark options={solid,fill=black}]
table [row sep=\\]{%
0.154303349962092	0.0170916215456275 \\
0.0816496580927726	0.010342901970135 \\
0.0399043442233811	0.00727127964708976 \\
0.0197027601559775	0.0043151127858851 \\
} node [below right, pos=0] {$2$};
\addplot [semithick, dashed]
table [row sep=\\]{%
0.154303349962092	0.38 \\
0.0197027601559775	0.048521622253248\\
};
\addplot [semithick, blue, mark=triangle, mark size=3, mark options={solid,fill=black}, forget plot]
table [row sep=\\]{%
0.154303349962092	0.425240852583659 \\
0.0816496580927726	0.221657173185071 \\
0.0399043442233811	0.111496071706317 \\
0.0197027601559775	0.055146233736575 \\
} node [below right, pos=0] {$1$};
\addplot [semithick, blue, mark=triangle, mark size=3, mark options={solid,fill=black}, forget plot]
table [row sep=\\]{%
0.154303349962092	0.48000984522911 \\
0.0816496580927726	0.248452507982595 \\
0.0399043442233811	0.11976428151818 \\
0.0197027601559775	0.0564551595222577 \\
} node [above right, pos=0] {$2$};
\end{axis}
\end{tikzpicture}
  \caption{Convergence rates for the two fractures in Example 2. The error is measured
    as the relative errors, $\varepsilon_{\Gamma_i^+}([\vec u_h], [\vec u])$ and
    $\varepsilon_{\Gamma_i^+}(\vec\lambda_h, \vec\lambda)$, for each fracture,
    $\Gamma_i^+$, as given in Equation~\ref{eq:relative_error}. The numbering of the
    lines correspond to the fracture number. \label{fig:ex2_rate}}
\end{figure}

The errors are calculated by comparison to a reference solution that has 500 thousand
degrees of freedom. The relative errors $\varepsilon_{\Gamma_i^+}([\vec u_h], [\vec u])$
and $\varepsilon_{\Gamma_i^+}(\vec\lambda_h, \vec\lambda)$ for the two fractures,
$i=1, 2$, are shown in Figure~\ref{fig:ex2_rate}. We observe first-order convergence for
the displacement jump, while the Lagrange multiplier shows a somewhat reduced order of
convergence. Finally, Table~\ref{table:newton_count} shows the number of Newton
iterations for each mesh, and we do not observe any significant increase in the number
of Newton iterations as the mesh is refined.

\begin{table}[ht]
  \centering
  \caption{The number of Newton iterations used for the different mesh sizes in Example
    1 and 2. The number of fracture faces and number of cells are given in the
    table.\label{table:newton_count}}
  \vspace{2mm}
  \label{tab:lam1}
  \scalebox{0.8}{
    \begin{tabular}{rrcrrc}
      \toprule
      \multicolumn{3}{c}{Example 1} & \multicolumn{3}{c}{Example 2}\\%
      \cmidrule(lr){1-3} \cmidrule(lr){4-6}
      \# fracture faces & \# cells & \# iterations & \# fracture faces & \# cells & \# iterations\\ 
      \midrule
      176\hspace{20pt} & 4538 & 5 & 126\hspace{20pt} & 449 & 3\\ 
      346\hspace{20pt} & 17482 & 4 & 450\hspace{20pt} & 1878 & 3\\ 
      682\hspace{20pt} & 68330 & 5 & 1884\hspace{20pt} & 11825 & 3\\ 
      1356\hspace{20pt} & 270756 & 7 & 7728\hspace{20pt} & 160417& 4\\ 
      \bottomrule
    \end{tabular}
  }
\end{table}

\subsection{Example 3}
In this example, we consider the same domain and material parameters as in Example 1,
but add a fluid. The permeability of the fluid is
$\mathcal K = 10^{-8}$~m$^2$Pa$^{-1}$s$^{-1}$, the storage coefficient is
$c_0=1\cdot10^{-10}$~Pa$^{-1}$, and the Biot coefficient is $\alpha=1$. The initial
displacement and pressure is set to zero, and the end time of the simulation is set to
$T = 5 c_0 H^2/ \mathcal K$, where $H=1$ m is the height of the domain. For the fluid,
we enforce homogeneous Neumann conditions on all sides except the left boundary, where a
zero pressure condition is given. For the mechanics, the left and the right boundaries
are given a homogeneous Neumann condition, and the bottom boundary is given a zero
Dirichlet condition. The top boundary is given a time varying boundary condition given
by
\begin{equation*}
  \vec g_{u,D}(\vec x, t) =
  \begin{cases}
    [0.005\text{ m}, -0.002 \text{ m}]^\top 2 t / T, & t < T / 2\\
    [0.005\text{ m}, -0.002 \text{ m}]^\top, & t \ge T
  \end{cases}
  \quad \vec x~\text{on top boundary.}
\end{equation*}
This condition enforces a linear increase of the boundary condition values in the first
half of the simulation, and after the boundary condition reaches the same value as in
Example 1 we keep it constant for the remainder of the simulation.

\begin{figure}[h]
  \centering \def\svgwidth{\textwidth}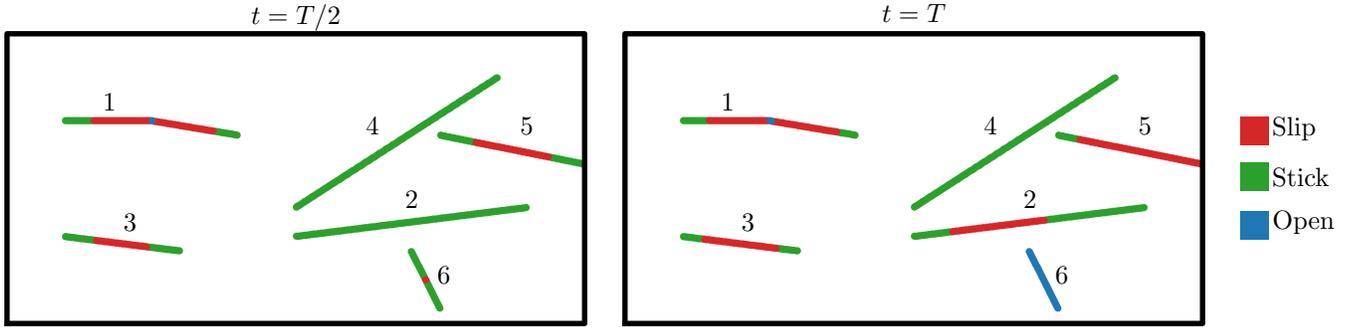
  \caption{The fractured domain from Example 3. The black box corresponds to the domain
    boundaries, while the fractures, labeled by a number, are represented by lines.  The
    fractures are colored in three colors; segments that slip are red, segments that
    stick are green, and segments that open are blue. The left figure shows the
    deformation half-way through the simulation, while the right figure, shows the
    deformation at the end time of the simulation. \label{fig:biot_slip_stick}}
\end{figure}

In the first half of the simulation, the domain is compressed and the fluid pressure in
the domain increases. In Figure~\ref{fig:biot_slip_stick}, we plot the regions of the
fractures that slip, stick and open for the two times $t=T/2$ and $t = T$. At the time
$t = T/2$, the boundary condition for the mechanics is the same as in Example 1,
however, we see considerable differences in the fracture displacement (compared to
Figure~\ref{fig:slip_stick}) that are due to the introduction of the fluid. After
$t=T/2$ the displacement boundary condition is fixed at $[0.005, -0.002]^\top$ m, and
the fluid pressure decrease due to the zero pressure condition on the left
boundary. This causes further deformation of the fractures, and at the end of the
simulation the pressure in the whole domain is relatively close to zero, and the
solution is approximately equal the solution in Example 1.
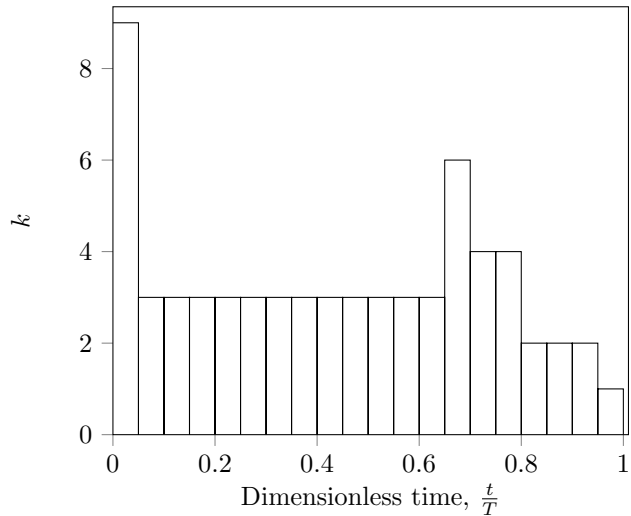
\begin{figure}
  \centering \begin{tikzpicture}

\begin{axis}[
tick align=outside,
tick pos=left,
x grid style={white!69.01960784313725!black},
xlabel={Dimensionless time, $\frac{t}{T}$},
xmin=0, xmax=1.01,
y grid style={white!69.01960784313725!black},
ylabel={$k$},
ymin=0, ymax=9.35
]
\addplot [const plot, draw=black]
table [row sep=\\,x expr=\thisrowno{0}-0.05]{%
0.05	9 \\
0.1	0 \\
0.1	3 \\
0.15	0 \\
0.15	3 \\
0.2	0 \\
0.2	3 \\
0.25	0 \\
0.25	3 \\
0.3	0 \\
0.3	3 \\
0.35	0 \\
0.35	3 \\
0.4	0 \\
0.4	3 \\
0.45	0 \\
0.45	3 \\
0.5	0 \\
0.5	3 \\
0.55	0 \\
0.55	3 \\
0.6	0 \\
0.6	3 \\
0.65	0 \\
0.65	3 \\
0.7	0 \\
0.7	6 \\
0.75	0 \\
0.75	4 \\
0.8	0 \\
0.8	4 \\
0.85	0 \\
0.85	2 \\
0.9	0 \\
0.9	2 \\
0.95	0 \\
0.95	2 \\
1	0 \\
1	1 \\
1.05	0 \\
};
\end{axis}

\end{tikzpicture}
  \caption{The number of Newton iterations at each time step for the simulation in
    Example 3.\label{fig:ex3_newton_it}}
\end{figure}

The number of Newton iterations for each time step is shown in
Figure~\ref{fig:ex3_newton_it}. For most iterations, the Newton solver converges in
three iterations. The increase in the number of iterations needed at time $t/T = 0.65$
is believed to be related to Fracture~6 changing behavior from sliding to opening. At
the end of the simulation the system is close to steady state, and the Newton solver
converges in one iteration.

\subsection{Example 4}
In this example, the same setup as in Example 2 is used, but a fluid is included. The
domain is sealed for the fluid, i.e., homogeneous Neumann conditions, for all sides
except the top boundary, which is given a Dirichlet condition $g_{p,D} = 0$~Pa.  The
permeability is $\mathcal K = 10^{-8}$~m$^2$Pa$^{-1}$s$^{-1}$, the storage coefficient
$c_0=1\cdot10^{-10}$~Pa$^{-1}$, and the Biot coefficient $\alpha=1$. The initial
displacement and pressure is set to zero.

Without the fractures, this setup is equivalent to a consolidation problem, which can be
found in standard textbooks~\cite{jaeger2007fundamentals}. When the load is applied to
the top surface at time $t=0$, there is an instantaneous increase in the pore pressure
in the domain. The fluid will then drain slowly out from the top surface and finally
relax back to the initial condition. As this process occurs, the domain will continue to
deform vertically increasing the mechanical load on the fractures, which causes them to
slip.  Twenty time steps are taken, and the simulation is stopped after $625$ minutes,
at which time, for practical purposes, equilibrium is reached.

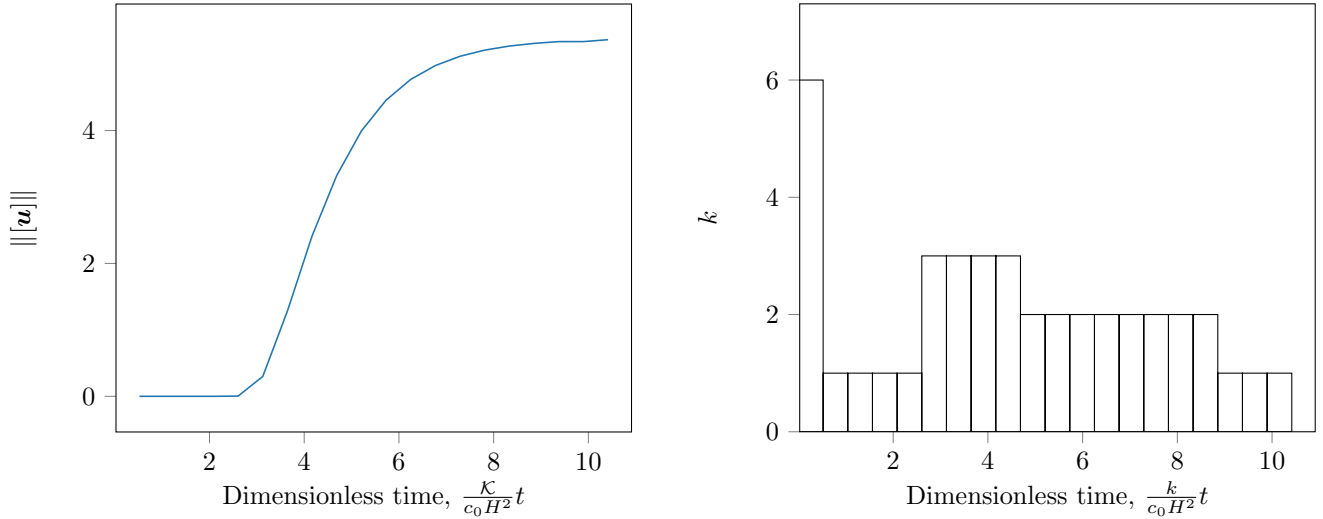
\begin{figure}
  \centering \begin{tikzpicture}

\definecolor{color0}{rgb}{0.12156862745098,0.466666666666667,0.705882352941177}

\begin{axis}[
tick align=outside,
tick pos=left,
x grid style={white!69.01960784313725!black},
xlabel={Dimensionless time, $\frac{\mathcal K}{c_0H^2}t$},
xmin=0.0260416666666667, xmax=10.9114583333333,
y grid style={white!69.01960784313725!black},
ylabel={$\lVert [\vec u] \rVert$},
]
\addplot [semithick, color0, forget plot]
table [row sep=\\]{%
0.520833333333333	2.15589487693912e-08 \\
1.04166666666667	2.1559043439312e-08 \\
1.5625	2.15591075933605e-08 \\
2.08333333333333	2.15591611159228e-08 \\
2.60416666666667	0.00326798255303682 \\
3.125	0.298558768068344 \\
3.64583333333333	1.28349309743603 \\
4.16666666666667	2.4141187814872 \\
4.6875	3.3239856069704 \\
5.20833333333333	3.99456689908579 \\
5.72916666666667	4.4574491505635 \\
6.25	4.76964942615286 \\
6.77083333333333	4.97862041131037 \\
7.29166666666667	5.11808665907887 \\
7.8125	5.21061683465126 \\
8.33333333333333	5.27170188341438 \\
8.85416666666667	5.31203911065458 \\
9.375	5.33879230248184 \\
9.89583333333333	5.33880584454838 \\
10.4166666666667	5.36764928501574 \\
};
\end{axis}

\end{tikzpicture} \hfill
  \begin{tikzpicture}

\definecolor{color0}{rgb}{0.12156862745098,0.466666666666667,0.705882352941177}
\definecolor{color1}{rgb}{1,0.498039215686275,0.0549019607843137}

\begin{axis}[
tick align=outside,
tick pos=left,
x grid style={white!69.01960784313725!black},
xlabel={Dimensionless time, $\frac{k}{c_0H^2}t$},
xmin=0.0260416666666667, xmax=10.9114583333333,
y grid style={white!69.01960784313725!black},
ylabel={$k$},
ymin=0, ymax=7.3
]
\addplot [const plot, draw=black]
table [row sep=\\,x expr=\thisrowno{0}-0.520833333333333]{%
0.520833333333333	6 \\
1.04166666666667	0 \\
1.04166666666667	1 \\
1.5625	0 \\
1.5625	1 \\
2.08333333333333	0 \\
2.08333333333333	1 \\
2.60416666666667	0 \\
2.60416666666667	1 \\
3.125	0 \\
3.125	3 \\
3.64583333333333	0 \\
3.64583333333333	3 \\
4.16666666666667	0 \\
4.16666666666667	3 \\
4.6875	0 \\
4.6875	3 \\
5.20833333333333	0 \\
5.20833333333333	2 \\
5.72916666666667	0 \\
5.72916666666667	2 \\
6.25	0 \\
6.25	2 \\
6.77083333333333	0 \\
6.77083333333333	2 \\
7.29166666666667	0 \\
7.29166666666667	2 \\
7.8125	0 \\
7.8125	2 \\
8.33333333333333	0 \\
8.33333333333333	2 \\
8.85416666666667	0 \\
8.85416666666667	2 \\
9.375	0 \\
9.375	1 \\
9.89583333333333	0 \\
9.89583333333333	1 \\
10.4166666666667	0 \\
10.4166666666667	1 \\
10.9375000000000	0 \\
}\closedcycle;
\end{axis}

\end{tikzpicture}
  \caption{Results from Example 4. The length scale $H=600$~m is the height of the
    domain. Left: Slip distance as a function of time. Right: Convergence of the Newton
    solver. Each time step is represented by a rectangle, and the $y$-axis shows the
    number of Newton iterations needed until the convergence tolerance is
    reached.\label{fig:ex4_slip}}
\end{figure}

The slip over time is plotted in Figure~\ref{fig:ex4_slip}. Initially, the pore pressure
carries most of the applied load, and the fractures are not sliding. As the fluid drains
and the domain deforms, the tangential part of the Lagrange multiplier on the fractures
increases, and after approximately 150 minutes, the fractures start to slide.  The
sliding then gradually slows down and qualitatively reaches the solution of the drained
medium, i.e., the solution from Example 2. There are small differences between the
solution from this example at the final time and the solution of Example 2, which are
caused by the use of a dynamic friction model in this example and a static friction
model in Example 2.

The number of iterations needed for convergence of the Newton solver at each time step
is shown in Figure~\ref{fig:ex4_slip}. For the first time step, 6 Newton iterations are
needed, which is twice as many as for any of the other time steps. It is well known that
the Newton strategy is very sensitive to the initial guess. A naive choice generally
results in an increase in the required number of Newton iterations for smaller mesh
sizes. However, either in a dynamic or a multilevel context, there are good options to
set the initial guess~\cite{hueber2008primal-dual,wohlmuth2011variationally}. In this
case, the initial condition at $t=0$ is ($\vec\lambda=\vec 0$ and $\vec u=\vec 0$),
which assigns all subfaces to the noncontact set, $\mathcal I_n$, while those at the
first time step belong to the sticking set $\mathcal I_\tau$ (see
Equation~\eqref{eq:sliding_sets}). During the dynamic sliding, the initial guess (the
solution from the previous time step) gives a good approximation of the solution in the
current time step, and thus, fewer iterations are needed. As the fractures start to
slide at time step six, a few Newton iterations are needed for convergence. However,
when approaching steady state, the algorithm predicts the correct slip in just one
iteration.

\section{Conclusion}
In this paper, we present an approach for solving the poroelastic Biot equations in a
fractured domain. A classical hybrid formulation for contact mechanics is combined with
a finite volume discretization for poroelasticity. The fractures are modeled as internal
contact boundaries and are governed by a nonpenetration condition in the normal
direction and a Coulomb friction law in the tangential direction. The inequalities in
the contact conditions are handled by a semismooth Newton method. The finite volume
discretization has several advantages for these types of problems. The cell-center
colocation of the discrete displacement and pressure variables gives a sparse linear
system, efficient data structures, and no need for staggered grids. Moreover, the
contact conditions are obtained naturally in the discretization as a condition per
subface in the local systems. Thus, these conditions can be treated in an equivalent
manner to boundary conditions on the external boundary.  Finally, there is no need for
special treatment of the contact conditions in the poroelastic case versus the purely
elastic case, as the correct pressure contribution to the effective stress is obtained
in the local system.

We showed that the hybrid formulation coupled with the finite volume discretization
handles a given spatially varying coefficient of friction. The formulation is also
suitable for other friction models such as rate and state friction or
temperature-dependent coefficient of friction.

Four numerical examples illustrate the method's robustness and applicability to
difficult cases. By comparison to a reference solution, the discrete solution shows
first-order convergence in displacements and slightly less than first-order convergence
for the Lagrange multipliers. We also show that the method handles singularity in the
solution resulting from a piecewise linear fracture with a kink. Finally, a 3d example
is presented where we study the effect of the fluid pressure on the solution.

The model presented in this work is limited to fluid flow in the matrix. A natural
extension is to include fluid flow also in the fractures. The fluid pressure in the
fracture will then act as a force on the fracture sides, effectively reducing the normal
traction. Experiments have also shown that asperities along fracture surfaces can have a
very important effect on both the opening and sliding of fractures. These effects can be
included by adding a nonlinear deformation model to the fractures. The advantage of our
framework is that any nonlinear extensions to the model can be included in the same
Newton iteration, which might be crucial for the convergence of the resulting scheme.

\printbibliography%

\appendixtitleon
\begin{appendices}
\section{}\label{sec:appendix}
\begin{figure}
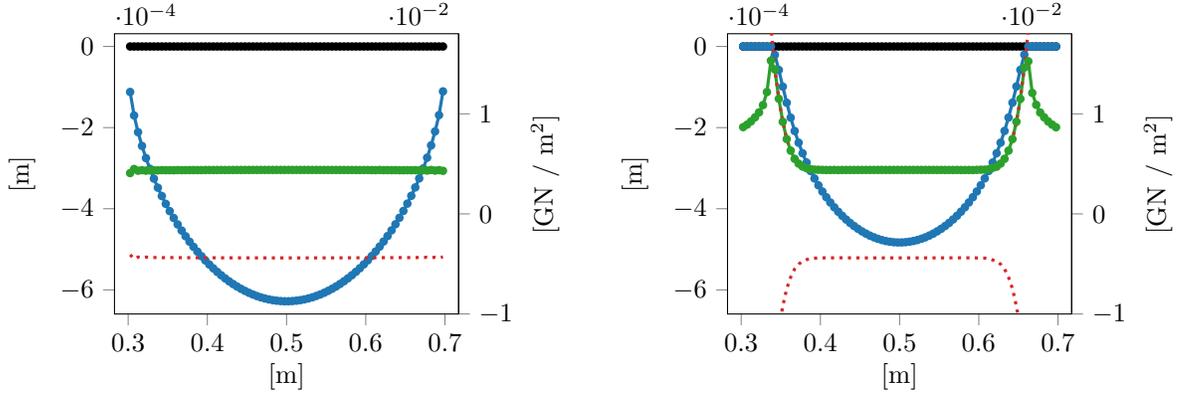

  \centering \def\figwidth{0.35\textwidth}
  \begin{subfigure}{\figwidth}
    \def\svgwidth{\textwidth} \import{fig/example0/}{traction_no_reg0.tex}
  \end{subfigure}
  \hspace{0.1\textwidth}
  \begin{subfigure}{\figwidth}
    \def\svgwidth{\textwidth} \import{fig/example0/}{traction_reg0.tex}
  \end{subfigure}
  \caption{The shear component $\vec\lambda_\tau$ (green), friction bound
    $\pm F\lambda_n$ (red dashes), tangential displacement jump $[\vec u]_\tau$ (blue),
    and normal displacement jump $[\vec u]_n$ (black) for the fracture.  The dots
    correspond to the face-centered values. Left: Constant friction coefficient.  Right:
    Regularized coefficient.\label{fig:ex0_traction}}
\end{figure}

\begin{figure}
  \centering\def\figwidth{0.35\textwidth}
  \begin{subfigure}{\figwidth}
    \def\svgwidth{\textwidth} \begin{tikzpicture}
  \definecolor{blue}{rgb}{0.12156862745098,0.466666666666667,0.705882352941177}
  \definecolor{green}{rgb}{0.172549019607843,0.627450980392157,0.172549019607843}

\begin{axis}[
width=\svgwidth,
tick align=outside,
tick pos=left,
legend cell align={left},
legend entries={{$\pmb\lambda_h$},{$\vec u_h$},{1st order}},
legend style={at={(0.03,0.97)}, anchor=north west, draw=white!80.0!black},
x grid style={white!69.01960784313725!black},
xlabel={$h = 1 / (\#${number of mortar cells}$)$},
xmin=0.00563281539131769, xmax=0.0554784736033923,
xmode=log,
y grid style={white!69.01960784313725!black},
ylabel={Relative error},
ymin=5e-3, ymax=1e-1,
ymode=log
]
\addlegendimage{mark=square, green}
\addlegendimage{mark=triangle, blue}
\addlegendimage{thick, dashed}
\addplot [thick, dashed]
table [row sep=\\]{%
0.05	0.081 \\
0.00625	0.010125 \\
};
\addplot [semithick, blue, mark=triangle, mark size=3, mark options={solid,fill=black}]
table [row sep=\\]{%
0.05	0.0688071161328575 \\
0.025	0.0335930707935884 \\
0.0125	0.0157605044410626 \\
0.00625	0.00675523480898647 \\
};

\addplot [semithick, green, mark=square, mark size=3, mark options={solid,fill=black}]
table [row sep=\\]{%
0.05	0.00795727676899801 \\
0.025	0.0103901841489542 \\
0.0125	0.0110530066447302 \\
0.00625	0.0071752439085591 \\
};
\end{axis}

\end{tikzpicture}
  \end{subfigure}
  \hspace{0.1\textwidth}
  \begin{subfigure}{\figwidth}
    \def\svgwidth{\textwidth} \begin{tikzpicture}
  \definecolor{blue}{rgb}{0.12156862745098,0.466666666666667,0.705882352941177}
  \definecolor{green}{rgb}{0.172549019607843,0.627450980392157,0.172549019607843}

\begin{axis}[
width=\svgwidth,
tick align=outside,
tick pos=left,
x grid style={white!69.01960784313725!black},
xlabel={$h = 1 / (\#${number of mortar cells}$)$},
xmin=0.00563281539131769, xmax=0.0554784736033923,
xmode=log,
y grid style={white!69.01960784313725!black},
ylabel={Relative error},
ymode=log
]
\addplot [thick, dashed]
table [row sep=\\]{%
0.05	0.081 \\
0.00625	0.010125 \\
};
\addplot [semithick, blue, mark=triangle, mark size=3, mark options={solid,fill=black}]
table [row sep=\\]{%
0.05	0.0637451962957316 \\
0.025	0.0302032787874455 \\
0.0125	0.0144363450147929 \\
0.00625	0.00646340210531151 \\
};
\addplot [semithick, green, mark=square, mark size=3, mark options={solid,fill=black}]
table [row sep=\\]{%
0.05	0.10221229712508 \\
0.025	0.058798459692563 \\
0.0125	0.0310782136974664 \\
0.00625	0.0122787171031989 \\
};
\end{axis}

\end{tikzpicture}
  \end{subfigure}
  \caption{Convergence rates for the jump $[\vec u_h]$ and Lagrange multiplier
    $\vec\lambda_h$ on $\Gamma^+$. The y-axis shows the relative error
    $\lVert \vec v_h - \vec v\rVert_{\Gamma^+} / \lVert \vec v\rVert_{\Gamma^+},\ 
    \vec
    v\in \{[\vec u],\vec\lambda\}$. Left: The convergence of the unregularized
    solution. Right: The convergence of the regularized solution. 
    \label{fig:ex0_rate}}
\end{figure}
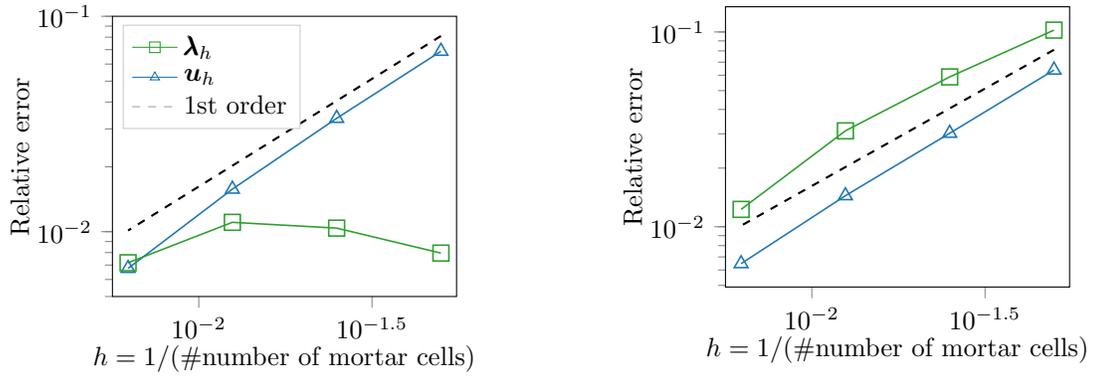

When a fracture slides or opens, the linear elastic stress will contain a singularity at
the fracture tips~\cite{lazzarin1996unified}, which causes challenges for any numerical
method.  We illustrate this in Figure~\ref{fig:ex0_traction}, where we plot the typical
stress and displacement profiles for a sliding fracture and a constant friction
coefficient $F=0.5$. We observe small oscillations in the Lagrange multiplier around the
tips of the fracture. The issue is that as we approach the fracture tips, an
infinitesimal change in the displacement jump will induce an infinite change in the
stress. These oscillations are reflected in the errors plotted in
Figure~\ref{fig:ex0_rate}, where the error rate for the Lagrange multiplier
deteriorates. Convergence is not seen in the Lagrange multiplier.  Because the face
traction values away from the fracture tips are almost constant, the error in this
region is very small, and thus, the error in the Lagrange multiplier is completely
dominated by the oscillations near the tips. Note that the convergence rates for the
displacement jump is of order 1, as expected. To study the convergence of the Lagrange
multiplier, we can regularize the solution by increasing the friction bound smoothly in
a small region around the tips.  In this example, this is done by setting
\begin{equation*}
  F(\vec x) = 0.5 (1 + 10 \exp(-800\text{ m}^{-2}D(\vec x)^2)) \qquad \vec x\in 
  \Gamma^+,
\end{equation*}
where $D(\vec x)$ is the distance from $\vec x$ to the tips of the fracture. As seen in
Figure~\ref{fig:ex0_traction}, this arrests the fracture before the tip, and the added
regularity gives first-order convergence in both the Lagrange multiplier and
displacements, as shown in Figure~\ref{fig:ex0_rate}.

\begin{figure}
  \centering\def\figwidth{0.35\textwidth}
  \begin{subfigure}[t]{\figwidth}
    \includegraphics[width=\textwidth]{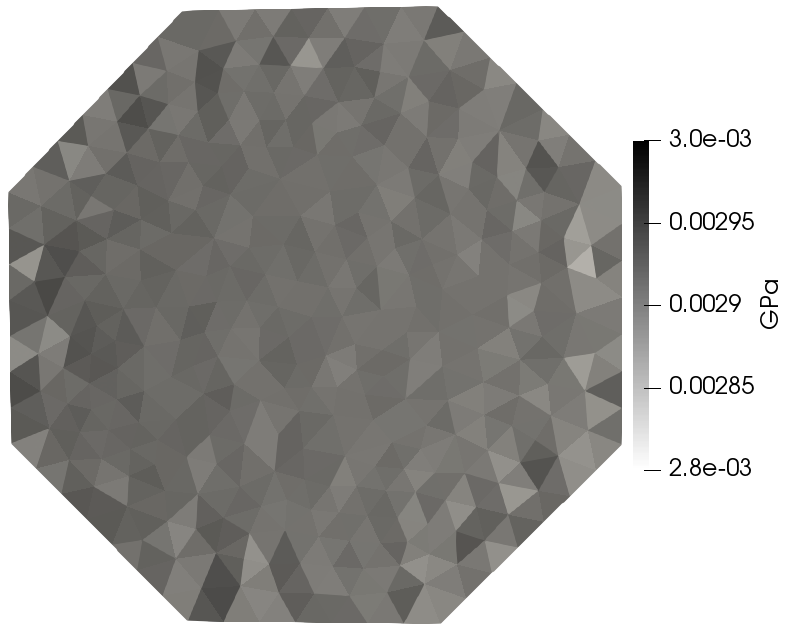}
  \end{subfigure}
  \hspace{0.1\textwidth}
  \begin{subfigure}[t]{\figwidth}
    \def\svgwidth{\textwidth} \import{fig/example0/}{oscillations.tex}
  \end{subfigure}
  \caption{Oscillations in the normal component of the Lagrange multiplier. Left: The
    negative normal component of the Lagrange multiplier $-\vec\lambda$ on the
    fracture. Right: The negative normal component of the Lagrange multiplier
    $-\vec\lambda$, where the $x$-axis is the radial distance from the fracture center,
    i.e., the center is at $x=0$, while the tip is at
    $x\approx150$~m.\label{fig:oscillations_3d}.}
\end{figure}

The worst oscillations that we have encountered in 3d using our finite volume scheme
coupled with the hybrid formulation are shown in Figure~\ref{fig:oscillations_3d}.  The
setup in this example is the same as the setup in Section~\ref{sec:example2} but with
only Fracture 1 and a constant coefficient of friction, $F=0.5$. Thus, we have sliding
reaching the tip of the fractures. The oscillations have an amplitude of approximately 5
percent from the mean traction and grow larger as we approach the fracture tips. As in
the 2d case, the displacement jump $[\vec u_h]$ is not effected significantly by these
oscillations.

Note that the singularity at the fracture tips is a challenge for any numerical method.
Similar oscillations for first- and second-order Galerkin finite elements are reported,
for example, by Garipov et al~\cite{Garipov2016discrete}$^\text{, Fig. 8}$, for a setup
where they study a single sliding fracture.
\end{appendices}
	
\end{document}